\documentclass[a4paper, 11pt, twoside, leqno]{amsart}

\usepackage[T1]{fontenc}
\usepackage[utf8x]{inputenc}
\usepackage[english]{babel}
\usepackage{amsmath, amsthm, amssymb, stmaryrd}
\usepackage{mathrsfs, esint, bbm}   % for special math fonts and symbols
\usepackage{fancyhdr}               % for headings
\usepackage{enumitem}               % for lists
\usepackage[all,cmtip]{xy}          % for commutative diagrams
\usepackage{graphicx}
\newtheorem{theorem}{Theorem}            % Bold title, italic text

\newtheorem{lemma}[theorem]{Lemma}
\newtheorem{prop}[theorem]{Proposition}

\theoremstyle{definition}              % Bold title, roman text
\newtheorem{definition}{Definition}

\theoremstyle{remark}                  % Italic title, roman text
\newtheorem{step}{Step}

\newtheorem{remark}{Remark}

\newcommand{\N}{\mathbb{N}}                 % natural numbers
\newcommand{\R}{\mathbb{R}}                 % real numbers
\newcommand{\C}{\mathbb{C}}                 % complex numbers
\renewcommand{\S}{\mathbb{S}}               % \S is now reserved for bbold S (e.g., spheres).

\let\oldv\v                                 % \oldv is the standard latex command \v. 
\renewcommand{\v}{\mathbf{v}}               % \v is now reserved for bold v and bold a.

\newcommand{\tOmega}{\tilde{\Omega}}
\renewcommand{\th}{\tilde{h}}
\renewcommand{\a}{\mathbf{a}}

\newcommand{\db}{\mathbf{d}}
\renewcommand{\d}{\mathrm{d}}               % Roman letters

\newcommand{\epsi}{\varepsilon}

\newcommand{\eps}{\epsi}
\newcommand{\argmin}{\hbox{argmin}}

\DeclareMathOperator{\sign}{sign}                                   % sign
\DeclareMathOperator{\curl}{curl} 
\DeclareMathOperator{\Argmin}{argmin}                                   % curl

\renewcommand{\div}{\mathrm{div}\,}      
\newcommand{\abs}[1]{\left| #1 \right|}                             % absolute value
\newcommand{\norm}[1]{\left\| #1 \right\|}                          % norm
\usepackage{color}
\usepackage[normalem]{ulem} % "\sout" for striking out in text mode
\numberwithin{equation}{section}
\numberwithin{definition}{section}
\numberwithin{theorem}{section}
\numberwithin{remark}{section}

\begin{document}

\title[Dimensional Reduction and defects for nematic liquid crystals]
{Dimensional Reduction and emergence of defects\\ in the Oseen-Frank model for nematic liquid crystals}
\author{Giacomo Canevari}
\address{Dipartimento di Informatica, 
  Universit\`a di Verona,
  Strada le Grazie 15, \mbox{37134} Verona, Italy.}
\email[G. Canevari]{giacomo.canevari@univr.it}
\author{Antonio Segatti}
\address{Dipartimento di Matematica ``F. Casorati'',
  Universit\`a di Pavia,
  Via Ferrata 1, \mbox{27100} Pavia, Italy.}
\email[A. Segatti]{antonio.segatti@unipv.it}

\date{\today}

%\dedicatory{\large{Dedicated to Pierluigi Colli on the occasion of his 65th anniversary, \\ with friendship and admiration}}

\begin{abstract}
 In this paper we discuss the behavior of the Oseen-Frank model for nematic liquid crystals in the limit of vanishing thickness. 
%  In this paper we relate the emergence of defect points in the two-dimensional Oseen-Frank model for nematic liquid crystals with non-zero degree Dirichlet boundary conditions, with the vanishing thickness limit of the Oseen-Frank model in three dimensions. 
 More precisely, in a thin slab~$\Omega\times (0,h)$ with~$\Omega\subset \R^2$ and $h>0$ we consider the one-constant approximation of the Oseen-Frank model for nematic liquid crystals. 
 We impose Dirichlet boundary conditions on the lateral boundary and weak anchoring conditions on the top and bottom faces of the
 cylinder~$\Omega\times (0,h)$. 
 The Dirichlet datum has the form $(g,0)$, where $g\colon\partial\Omega\to \mathbb{S}^1$ has non-zero winding number. 
 Under appropriate conditions on the scaling, in the limit as~$h\to 0$ we obtain a behavior that is similar to the one observed in the asymptotic analysis (see \cite{BBH}) of the two-dimensional Ginzburg-Landau functional. 
 More precisely, we rigorously prove the emergence of a finite number of defect points in $\Omega$ having topological charges that sum to the degree of the boundary datum. 
 Moreover, the position of these points is governed by a Renormalized Energy, as in the seminal results of Bethuel, Brezis and H\'elein~\cite{BBH}.
 
 \smallskip
 \noindent
 \medskip {\bf Keywords:} $\Gamma$-Convergence; Nematics; Dimension Reduction; Defects; Renormalized Energy.   
 
 \smallskip
 \noindent
 \medskip {\bf MSC:} 49J45 (82D30; 35Q56).

\end{abstract}

\maketitle
\tableofcontents

\section{Introduction}

We let $\Omega$ be a bounded, regular, simply connected domain in $\mathbb{R}^2$ and set $Q:=\Omega\times (0,1)$.
For any $\eps>0$ we consider the energy 
\begin{equation}
\label{eq:energyQ}
F_\eps(U) := \frac{1}{2}\int_{Q}\abs{\nabla_\eps U}^2\d x 
+ \frac{1}{2\eps^2}\int_{\Omega\times \left\{0,1\right\}}\abs{(\nu,U)}^2\d x',
\end{equation}
where $Q\ni x = (x',x_3)$, $\nu= \pm \hat{e}_3$ is the outer normal 
to~$\Omega\times \left\{0\right\}$ and~$\Omega\times\left\{1\right\}$, respectively, and the operator~$\nabla_\eps$ is defined as
\[
 \nabla_\eps:=\left(\frac{\partial}{\partial x_1},
 \, \frac{\partial}{\partial x_2}, 
 \, \frac{1}{\eta(\eps)}\frac{\partial}{\partial x_3}\right) \! , 
\]
depending on a parameter~$\eta(\eps)$ which is strictly positive 
and satisfies~$\eta(\eps)\xrightarrow{\eps\to 0}0$.
The map $U\colon Q\to \R^3$ belongs to the set of the admissible configurations
\begin{equation}
\label{eq:funct-frame}
\mathcal{A}_{G}:= \left\{V:Q\to \R^3: V\in H^1(Q;\R^3), \,\, \abs{V}=1\hbox{ a.e. in } Q,\,\, V = G \hbox{ on } \partial \Omega\times (0,1) \right\} \! .
\end{equation}
%{\BBB\textbf{Per uniformare la notazione, ho sostituito~$\tilde{g}$ con~$G$.}}
We assume that the boundary datum
$G\in H^{1/2}(\partial\Omega\times (0,1);\R^3)$
takes the form
\begin{equation}
\label{eq:bddatum}
 G(x^\prime, \,  x_3) = (g(x^\prime), \, 0)
 \qquad \textrm{for } (x^\prime, x_3)\in \partial\Omega\times (0, \, 1),
\end{equation}
where~$g$ is a given map in $g\in H^{1/2}(\partial\Omega;\S^1)$
(that is, $g\in H^{1/2}(\partial\Omega; \R^2)$ and~$\abs{g}=1$
a.e.~on~$\partial\Omega$ with respect to the length measure).
In particular, the map~$g$
has a well-defined topological degree on~$\partial\Omega$,
denoted as~$d:=\deg(g,\partial\Omega)$
(see Section~\ref{ssec:jaco} below for details).
We are interested in the case~$d \neq 0$ (see Remark \ref{oss:degree} for a (brief) discussion on the $d=0$ case);
for simplicity, throughout the paper assume that~$d > 0$.
% (if~$\deg(u; \partial U) < 0$, we can reduce
% to the case~$d > 0$ by replacing~$g$ with~$\bar g := (g_1, \, -g_2)$).
The functional~$F_\eps$ %emerges after a non dimensionalization procedure of
is a non-dimensional form of the %(one constant approximation) 
Oseen-Frank energy for nematic liquid crystals 
on a thin cylindrical domain,
in the one-constant approximation. We refer to the next Section~\ref{sec:motivation} for the details. 

The main aim of this paper is to discuss the behavior of the energy~$F_\eps$ as $\eps\to 0$. %in terms of $\Gamma$-convergence. 
For any $\eps>0$, the direct method of the calculus of variations guarantees the existence of a minimizer $U_\eps$. 
When $\eps\to 0$ we have the following heuristic on the behavior of $U_\eps$. On the one hand, since $\eta(\eps)\xrightarrow{\eps \to 0}0$, 
%one expects that 
the energy promotes that $\partial_{x_3}U_\eps\approx 0$, namely
$U_\eps$ is independent of $x_3$. 
In this sense, we are discussing the {\itshape dimension reduction} of the Oseen-Frank energy --- that is, we obtain in the limit an energy for maps defined on the two dimensional domain $\Omega$.\footnote{To be precise, the domain of the limit functional consists of maps defined in $Q$ that do not depend on the variable $x_3\in (0,1)$. In this sense, we might think the limit problem is set in~$\Omega$. See the Braides's book \cite[Chapter 14]{braides-beginner} for a presentation of dimensional reduction via~$\Gamma$-convergence.} 
On the other hand, the second term in the energy (named weak anchoring in the 
parlance of liquid crystals, see~\cite{Virga94})
% promotes also the fact that
favors configurations that satisfy
$\abs{\left(U_\eps,\nu\right)}\approx 0$ almost everywhere on $\Omega$. Consequently, since $\abs{U_\eps}=1$,
the two-dimensional projection~$\Pi(U_\eps)$ of~$U_\eps$
onto the~$(x_1, \, x_2)$-plane must satisfy
$\abs{\Pi(U_\eps)}\approx 1$ almost everywhere on $\Omega$.
% ($\Pi:\R^3\to \R^2$ is the orthogonal projection on the $x_3=0$ plane, see below).
As a result, when the boundary datum %$G$ is chosen in such a way that
%$G(\cdot, 0)=g(\cdot)\colon\partial\Omega\to \S^1$
$g\colon\partial\Omega\to \S^1$ has winding number $d\neq 0$, 
% we expect that in the limit we reproduce the behavior of the two-dimensional Ginzburg-Landau energy discussed in \cite{BBH}, \cite{JerrardSoner-GL}, \cite{AlicandroPonsiglione} ......
we expect that topological singularities will emerge in the limit,
as observed in the asymptotic behavior of the two-dimensional
Ginzburg-Landau functional.

The Ginzburg-Landau functional, in the simplified form 
considered by Bethuel, Brezis and H\'elein~\cite{BBH},
reads
\begin{equation} \label{GL}
 GL_\eps(u) := \int_{\Omega}\left(\frac{1}{2}\abs{D u}^2 + \frac{1}{4\eps^2}\left(1 - \abs{u}^2\right)^2\right) \d x^\prime
 \qquad \textrm{for } u\colon \Omega\subset\R^2\to\C.
\end{equation}
(Here,~$D$ denotes the usual Euclidean gradient in~$\R^2$.)
Minimizers of~\eqref{GL} subject to a boundary condition
$u_{|\partial\Omega} = g\in H^{1/2}(\partial\Omega;\S^1)$
have an interesting behavior as~$\eps\to 0$,
especially when~$d := \deg(g, \partial\Omega)\neq 0$.
Indeed, if the degree is nonzero (and, say, positive),
the boundary datum cannot 
be extended to a smooth map~$\Omega\to\S^1$ 
--- in fact, not even to a map in~$H^1(\Omega; \S^1)$ ---
because of topological obstructions. As a consequence,
the energy of the minimizers diverge as~$\eps\to 0$.
Nevertheless, it is still possible to obtain compactness for minimizers.
Bethuel, Brezis and H\'elein~\cite{BBH} showed that, up to
extraction of a subsequence, minimizers converge to a map 
$u^*\colon\Omega\setminus\{a_1^*, \, \ldots, \, a^*_d\}\to\S^1$ 
that is smooth except for finitely many singularities 
at the points~$a^*_j$. Once the boundary datum and the
singular points~$a^*_j$ are known, the limit map~$u^*$ is
uniquely determined, in terms of a (singular) Poisson problem;
in the language of~\cite{BBH}, $u^*$ is the canonical harmonic map
associated with~$a^*_j$ and~$g$ (see Definition~\ref{def:chm} below). 
Moreover, the points~$a_j^*$ minimize a function, the Renormalized Energy
$W_g = W_g\left(a_1^{*}, a_2^{*},\ldots, a_d^{*}\right)$
(see Section~\ref{sse:reno} below),
which expresses the interaction energy between the
topological singularities.
Over the decades, the analysis of~\cite{BBH} has been
generalized and extended to several different contexts,
including higher-dimensional domains 
(see e.g.~\cite{LinRiviere, BethuelBrezisOrlandi}),
evolution problems~\cite{SS-GF, BethuelOrlandiSmets-Annals},
and $\Gamma$-convergence results~\cite{JerrardSoner-GL, ABO2, AlicandroPonsiglione}. Other variants of the functional have been studied as well (see e.g.~\cite{SS-book}).

In this paper, we obtain a similar characterization 
for minimizers of~\eqref{eq:energyQ}, in the limit as~$\eps\to 0$,
under the assumption that
\begin{equation} \label{eq:ipeta}
 \sqrt{2}\eta(\eps) \leq \eps \qquad
 \textrm{for any } \eps > 0.
\end{equation}
(Minimizers of~$F_\eps$ do exist; see 
Section~\ref{ssec:funct_sp_min} for details.)
Mathematically, this regime of parameters correspond to a
dimensional reduction limit~$h(\eps)\to 0$. Moreover,
configurations that depend on the~$x_3$-variable are
more heavily penalized by the energy than 
anti-plane configurations. (For a discussion of the
physical model, see Section~\ref{sec:motivation}.)
% {\BBB For technical reasons, we also assume that
% \begin{equation} \label{eq:ipeta2}
%  \eps\mapsto\frac{\eta(\eps)}{\eps} \quad 
%  \textrm{is monotonically nondecreasing.}
% \end{equation}
% }

\begin{theorem} \label{th:minimizers}
 Assume that the boundary datum takes the form~\eqref{eq:bddatum},
 for some~$g\in H^{1/2}(\partial\Omega;\S^1)$ 
 with~$d := \deg(g, \partial\Omega) > 0$, 
 and that~\eqref{eq:ipeta} holds.
 Let~$U_\eps^*$ be a minimizer of~$F_\eps$ in the class~$\mathcal{A}_g$.
 Then, there exist a (non-relabelled) subsequence, 
 distinct points~$a^*_1, \, \ldots, \, a^*_d$ in~$\Omega$
 and a map $u^*\colon\Omega\setminus\{a_1^*, \, \ldots, \, a^*_d\}\to\S^1$ 
 such that the following properties hold:
 \begin{enumerate}[label=(\roman*)]
  \item minimizers~$U^*_\eps$ converge 
  (strongly in~$L^q(Q)$ for any~$q<+\infty$)
  to the map~$U^*(x) := (u^*(x^\prime), 0)$;
  \item the map~$u^*$ belongs to~$W^{1,p}(\Omega;\S^1)$
  for any~$p\in(1, \, 2)$,
  is smooth in~$\Omega\setminus\{a^*_1, \, \ldots, \, a^*_d\}$,
  has a topological singularity of degree~$1$ at each point~$a^*_j$,
  and satisfies~$u^* = g$ on~$\partial\Omega$;
  \item $u^*$ is, in fact, the canonical harmonic vector field
  associated to the singular points~$a^*_1$, \ldots, $a^*_d$
  and to~$g$;
  \item the points~$a^*_1$,\ldots $a^*_d$ minimize the Renormalized Energy~$W_g$.
%  {\BBB and there holds
  %\begin{equation}
   %\label{eq:energy_exp_intro}
   %F_\eps\left(U^*_\eps\right) = d\pi \abs{\log\eps} + W_g\left(a_1^{*}, a_2^{*},\ldots, a_d^{*}\right) + d\gamma + \mathrm{o}_{\eps\to 0}(1)
  %\end{equation}
 % where~$\gamma\in\R$ is a constant, defined in~\eqref{eq:conv_core}.}
 \end{enumerate}
\end{theorem}
% {\BBB Personalmente, eliminerei dall'enunciato l'espansione asintotica dell'energia (la parte in blu), 
% e il riferimento alla core energy~$\gamma$, perch\'e non \`e
% proprio elegantissimo che~$\gamma$ dipenda dalla
% sottosuccesione di~$\eps$. Per\`o, se preferisci, la teniamo.}
% {\GGG Ho levato l'espansione. E' commentata nel .tex}

\begin{remark}[On the degree of the boundary condition]
\label{oss:degree}
 The assumption that~$d = \deg(g, \, \partial\Omega) > 0$
 is not restrictive. Indeed, given any~$U\in\mathcal{A}_G$,
 the map $\bar{U}\colon Q\to\S^2$ given component-wise by 
 $\bar{U} := (U_1, \, -U_2, \, U_3)$ has exactly
 the same energy as~$U$, but its trace 
 on~$\partial\Omega\times (0, \, 1)$ has degree~$-d$. 
 
 The case $d=0$ is simpler. The main reason is that in this case the space
 \[
 H^1_{g}(\Omega;\mathbb{S}^1):=\left\{v\in H^1(\Omega;\R^2): v=g \textrm{ on } \partial \Omega \,\,\hbox{ and } \abs{v}=1 
 \textrm{ in }\Omega\right\}
 %\neq \emptyset.
 \]
 is nonempty.
 With this observation, the analysis follows by a routine application of the arguments in \cite[Chapter 14]{braides-beginner}. More precisely, 
 when $d=0$, the energy $F_\eps$\, $\Gamma$-converges (in the topology of $L^2(Q;\R^3)$, for instance) to the energy 
 \[
 F(U) = 
 \begin{cases}
\displaystyle\frac{1}{2}\int_{Q} \abs{\nabla U}^2 \d x \qquad &\textrm{ if } U\in \mathcal{X}_G,\\
\displaystyle +\infty \qquad & \textrm{ otherwise in } L^2(Q;\R^3),
 \end{cases}
 \]
 where 
 \begin{equation}
 \label{eq:X}
 \mathcal{X}_G := \left\{V\in H^1(Q;\R^3): \partial_{x_3}V =0,  \, \abs{(V,\nu)}=0, \ %\textrm{ and } 
 \abs{V}=1 \,\textrm{ in } Q \textrm{ and } V = G \textrm{ on }\partial\Omega\times (0,1) \right\}.
 \end{equation}
 Then, the analogue of Theorem \ref{th:minimizers} is a consequence of the Fundamental Theorem of $\Gamma$-con\-ver\-gence. 
 Note also that the space $\mathcal{X}_G$ is isomorphic to $H^1_{g}(\Omega;\mathbb{S}^1)$.
 Therefore, the limit energy can be thought of as the two-dimensional Oseen Frank energy in the one-constant approximation.
\end{remark}

\begin{remark} \label{rk:nonunique}
 Theorem~\ref{th:minimizers} will follow as a corollary
 from the~$\Gamma$-convergence analysis contained
 in Section~\ref{sec:proofGamma} below. More precisely,
 we show that the sequence of functionals~$F_\eps$
 is compact in the sense of~$\Gamma$-convergence,
 and the limit functional(s) can be described in terms of
 the Renormalized Energy (see Theorem~\ref{th:Gamma}).
 We do \emph{not} identify a unique $\Gamma$-limit:
 it is possible that, along different subsequences~$\eps_k\to 0$
 and~$\eps_k^\prime\to 0$,
 the $\Gamma$-limits~$F := \lim_{k\to+\infty} F_{\eps_k}$ 
 and~$F^\prime := \lim_{k\to+\infty} F_{\eps_k^\prime}$ 
 might be different. However, the difference~$F - F^\prime$
 is (locally) constant, and minimizers of~$F$, $F^\prime$ 
 are the same. (See Sections~\ref{ssec:core} 
 and~\ref{ssec:Gamma} for details.)

 There are, however, instances in which we have a proper
 $\Gamma$-convergence result. 
 This is the case, for instance, when $\eta(\eps) = k\, \eps$
 for some constant~$k\in (0, \frac{1}{\sqrt{2}}]$.
 (Again, see Sections~\ref{ssec:core} 
 and~\ref{ssec:Gamma} for details.)
\end{remark}

In view of Theorem~\ref{th:minimizers}, in the regime
of parameters described by~\eqref{eq:ipeta},
the minimizers of~$F_\eps$ behave asymptotically as 
two-dimensional Ginzburg-Landau minimizers;
in particular, minimizers develop line singularities,
whose position is governed by the Ginzburg-Landau Renormalized Energy.
% Therefore, Theorem~\ref{th:minimizers}
% may be understood as a dimension reduction result.
For other dimension reduction results which apply
to liquid crystals, see e.g.~\cite{GMS15, Novack}.
(In these papers, the authors adopt a different modelling framework
for the liquid crystals, i.e. they work in the 
Landau-de Gennes theory as opposed to the Oseen-Frank theory.) 
% {\RRR Ci sarebbero anche Baumann \& Phillips, che hanno presentato 
% dei risultati (sempre riduzione dimensione LdG)
% ad una conferenza del~2021, ma non trovo un lavoro 
% scritto da nessuna parte. Posso provare a scrivere loro, se vuoi: 
% {\textbf{Ok, forse possiamo scrivere dopo aver sottomesso allegando il link ad arxiv. Poi, in caso, li citiamo nella revised version}}}

It would be interesting to extend this dimension reduction analysis to %the case of
the so-called Nematic Shells, namely rigid particles coated with a thin film of nematic liquid crystal. 
In this case, the simplest modeling approach that includes also the effects of the extrinsic geometry  would be to consider (see \cite{NapVer12E} and \cite{NapVer12L}) the following energy: 
\begin{equation}
\label{eq:shell}
E(v) = \frac{1}{2}\int_{\Sigma} \abs{D v}^2 + \abs{\mathfrak{B}v}^2 \d S.
\end{equation}
In the display above, $\Sigma$ is a two-dimensional closed (that is, compact, connected and without boundary) and oriented surface, isometrically embedded in $\R^3$; $v$ is a unit-norm tangent vector field on~$\Sigma$; $D$ is the covariant derivative on~$\Sigma$ and~$\mathfrak{B}$ is the shape operator.  
It is well known (see e.g. \cite{gamma-discreto}, \cite{AGM-VMO}, \cite{JerrardIgnat_full}) that, for closed and oriented surfaces in $\R^3$ with Euler Characteristic $\chi(\Sigma)\neq 0$, the incompatibility between the topology of $\Sigma$ and the constraints of unit norm and tangency is the trigger mechanism for the development of topological singularities that make the energy equal to $+\infty$. 
In this case, the minimization of \eqref{eq:shell} requires at first a relaxation of one of the constraints. 
Among the possibile strategies we recall the ``discrete-to-continuum approach'' in~\cite{gamma-discreto} and the Ginzburg-Landau approach of Ignat \& Jerrard  in~\cite{JerrardIgnat_full}. 
A third possibile strategy would be to consider, as in the present paper, a tubular neighborhood of $\Sigma$ and study the related dimensional reduction process via $\Gamma$-convergence. We expect that, under suitable assumptions on the scaling, one might obtain the results of Ignat \& Jerrard. 
We remark that in the case of a shell with %trivial topology (i.e. with 
$\chi(\Sigma)=0$ the energy \eqref{eq:shell} has been obtained in \cite{NapVer12E} exactly via a dimensional reduction procedure.\footnote{To be precise, the paper \cite{NapVer12E} deals with the full Oseen-Frank energy, not necessarily restricted to the one-constant approximation. Moreover, the arguments of~\cite{NapVer12E} 
does not use the theory of $\Gamma$-convergence. The asymptotic analysis as the thickness goes to zero is performed in the sense of pointwise convergence in~$\mathcal{X}_G$ (see \eqref{eq:X}). 
In the case of the energy $F_\eps$, this means 
% fixing $v\in \mathcal{X}_G$
% and then studying $\lim_{\eps\to 0}F_\eps(v)$
proving that~$F(V) = \lim_{\eps\to 0}F_\eps(V)$
for any~$V\in\mathcal{X}_G$
(see \cite[Proposition 2]{NapVer12E}). In the language of $\Gamma$-convergence, this would correspond to the existence of a (constant) recovery sequence. 
We refer to \cite{GMS17} for a dimensional reduction analysis via $\Gamma$-convergence for the Landau-de Gennes model on thin shells.}

Regarding the scaling assumption, it would be interesting to consider regimes of
parameters that do \emph{not} satisfy~\eqref{eq:ipeta}.
For instance, when~$\eta(\eps)$ is sufficiently large compared to~$\eps$,
the optimal configurations may depend on~$x_3$ in a substantial way
and they may exhibit other type of defects
--- for instance, surface defects, or `boojums', as they 
are known in the liquid crystals literature.
The emergence of surface defects in variational models
reminiscent of the Ginzburg-Landau functional has been
investigated, e.g., in~\cite{AlamaBronsardGolovaty}
in the context of liquid crystals and~\cite{IgnatKurzke}
in the context of micromagnetism.

\subsection{Notation} 

Given a smooth map $U:Q\to \R^3$, we write
\begin{equation}
\label{eq:ortogonal_decoU}
U = \sum_{i=1}^2 (U, \hat{e}_i) \hat{e}_i + (U,\hat{e}_3)\hat{e}_3 = \Pi(U) + U_{\perp}\hat{e}_3,
\end{equation}
where $\Pi:\mathbb{R}^3\to \mathbb{R}^2$ is the orthogonal projection on the $x_3=0$ plane. For brevity, we will often use a lowercase letter to indicate the projection of a vector field,
namely $u:=\Pi(U)$. 
We have
\begin{equation}
\label{eq:orto_deriv1}
\abs{\nabla_\eps U}^2 := \nabla_\eps U:\nabla_\eps U = \sum_{k=1}^2 \abs{\frac{\partial U}{\partial x_k}}^2 + \frac{1}{\eta^2(\eps)}\abs{\frac{\partial U}{\partial x_3}}^2,
\end{equation}
where, given two $(n\times m)$-matrices $A = a_{ij}$ and $B=b_{ij}$, we have denoted 
\[
A:B:= \sum_{i=i}^{n}\sum_{j=1}^{m}a_{ij}b_{ij},
\]
and, for $k=1,2$, 
\[
\frac{\partial U}{\partial x_k} = \frac{\partial u}{\partial x_k}+ \frac{\partial U_{\perp}}{\partial x_k}\hat{e}_3.
\] 
We denote by $D u$ the $(2\times 2)$-submatrix of the Jacobian matrix~$\nabla U$ whose columns are~$\frac{\partial u}{\partial x_k}$, for~$k=1,2$. 
Then, \eqref{eq:orto_deriv1} becomes 
\begin{equation}
\label{eq:orto_deriv}
\abs{\nabla_\eps U}^2 = \left(\abs{D u}^2 + \sum_{k=1}^2 \abs{\frac{\partial U_{\perp}}{\partial x_k}}^2\right) 
+  \frac{1}{\eta^2(\eps)}\abs{\frac{\partial U}{\partial x_3}}^2.
\end{equation}

Given a function $U\in H^1(Q;\R^3)$ and an open set $A\subset \Omega$,
we consider $B:=A\times (0,1)\subset Q$ and we set
\begin{equation}
 \label{eq:restr_energy}
 F_\eps\left(U; B\right):= \frac{1}{2}\int_{B}\abs{\nabla_\eps U}^2 \d x + \frac{1}{2\eps^2}\int_{A\times \left\{0,1\right\}}\abs{(\nu, U)}^2\d x'.
\end{equation}
 
%{\BBB Find a notation for the scalar product: for the moment we are using $(u,v)$ for the scalar product in $\R^2$ and for the scalar product in $\R^3$. The same symbol $(\cdot,\cdot)$ is used to denote the coordinates for the point $x\in Q$ with respect to the canonical basis in $\R^3$.
%\textbf{\`E davvero cos\`i problematico?} se non ti piacciono le tonde,
%possiamo usare~$\langle\cdot , \, \cdot\rangle$ per i prodotti scalari in~$\R^2$ e~$\R^3$, ma per il resto\ldots}
%
Unless otherwise stated, the notation $B_r(a)$
will denote a \emph{two-dimensional} disk,
of center~$a\in\R^2$ and radius~$r>0$. In case the center is the origin,
we write~$B_r := B_r(0)$. (Only in the proof of Lemma~\ref{lemma:nonempty} 
will we consider three-dimensional balls.)

\section{Nematic liquid crystals on a thin slab: the Oseen-Frank Energy}
\label{sec:motivation}

Liquid crystals, as the name suggests, are intermediate states of matter between solids and liquids.
liquid crystals are fluids, but they have anisotropic optical properties 
and their molecules can re-orient themselves when they are subject to an external electromagnetic field. For these reasons, liquid crystals are materials of choice for the multi-billion dollar industry of displays, and 
% their surprising ability of orienting  external electromagnetic fields have generated the multibillion dollar industry of 
liquid crystals displays %(LCDs) 
are still commonly found in personal computers, laptops, cellphones etc.

% There are at least three different types of liquid crystals: Nematic, Smectic and Cholesteric. 
In this paper we concentrate on a particular liquid crystal phase, the so-called nematic one --- which, incidentally, is the most commonly employed in the display industry. 
The name `nematic' derives from the greek $\nu \eta \mu \alpha$ (thread) and refers to a particular type of thread-like topological defects that this class of liquid crystals exhibits. 
In the nematic phase the (rod-like) molecules retain some orientational order but no positional order. 
As a consequence, one can describe the mean orientation of the molecules at a point $x$ (here we are neglecting all the possible time dependence) by a unit norm vector field. 

A well-established mathematical description of nematic liquid crystals is given by the Oseen-Frank variational theory (\cite{oseen33}, \cite{zocher33}, \cite{frank58}). 
We consider a thin film of nematic liquid crystals
contained in the volume~$\tilde{Q}= \tilde{\Omega}\times (0,\tilde{h})$,
where~$\tilde{\Omega}$ is a bounded, smooth, and simply connected domain in~$\R^2$ and $\tilde{h}>0$.
% As already observed, we describe the orientation of the optic axis by a map, named director, 
In the Oseen-Frank theory, the average orientation of the optic axis of the molecules is described by a unit vector field, i.e. a map 
\[
\tilde{U}:\tilde{Q}\to \mathbb{S}^2,
\]
where $\S^2:= \left\{x\in \R^3: \abs{x}=1\right\}$ is the unit sphere. 

The free energy density, $\sigma$, is assumed to 
be a function of~$\tilde{U}$ and of~$\nabla \tilde{U}$. 
More precisely, assuming that~$\sigma$ is a quadratic 
form of the gradient~$\nabla\tilde{U}$
and requiring it to be
{\itshape frame indifferent}, {\itshape even} 
(that is, invariant under the transformation $\tilde{U}\mapsto -\tilde{U}$) 
and {\itshape positive definite}, 
% and looking for a quadratic function of $\nabla\tilde{U}$,
Frank obtained the following expression for the free energy density:
\begin{equation}
\label{eq:frank_ed}
\begin{split}
 2\sigma\left(\tilde{U},\nabla \tilde{U}\right) 
 = k_1 (\div \tilde{U})^2 &+ k_2\left(\tilde{U}\cdot \curl\tilde{U} \right)^2 + k_3 \left(\tilde{U}\times \curl\tilde{U}\right)^2 \\
 &+ (k_2+k_4)\left(\textrm{tr}\left(\nabla \tilde{U}\right)^2-(\div\tilde{U})^2\right),
\end{split}
\end{equation}
Here, $k_1$, $k_2$, $k_3$, $k_4$ are elastic constants\footnote{The constants $k_1$, $k_2$, $k_3$, $k_4$ are of order magnitude of around $5\times 10^{-2}N$ (see, e.g., \cite{dg}).}, which
% are positive constants while $k_4\in \R$. The constants $k_i$ %$i=1\ldots,4$ 
which we assume to satisfy the so-called (strict) Ericksen inequalities~\cite{Eri66}: 
\begin{equation}
\label{eq:EI}
k_1 >0, \quad k_2>0,\quad k_3>0,\quad k_2>\abs{k_4}, \quad 2k_1>k_2+k_4.
\end{equation}
Equation~\eqref{eq:EI} is a necessary and sufficient
condition for the inequality
\begin{equation}
\label{eq:coerc}
\sigma(\tilde{U}, \, \nabla \tilde{U})\ge C\abs{\nabla\tilde{U}}^2
\end{equation}
to be satisfied, for some $C = C(k_1, k_2, k_3, k_4)>0$ (see e.g.~\cite{BallCime}). The constants $k_1$, $k_2$, $k_3$, called 
{\itshape{splay, twist, bend,}} respectively, 
account for the prototypical elastic deformations
that a nematic liquid crystal might experience
% weight the contribution to the free energy of the prototypical deformations that a nematic liquid crystal might experience
% splay, twist, bend 
(we refer to e.g.~\cite{Virga94} for a complete and detailed description). 
The fourth term 
% $(k_2 + k_4)\left(\textrm{tr}\left(\nabla \tilde{U}\right)^2-(\div\tilde{U})^2\right)$ 
is called {\itshape saddle-splay}.
The resulting Oseen-Frank free energy is thus the following:
\begin{equation}
\label{eq:OF_free}
W(\tilde{U}):= \int_{\tilde{Q}}\sigma (\tilde{U},\nabla\tilde{U}) \, \d \tilde{y}.
\end{equation}
The integration variable~$\tilde{y}\in \tilde{Q}$ 
will be written as $\tilde{y} = (\tilde{y}',\tilde{y}_3)$, 
with~$ \tilde{y}'=(\tilde{y}_1,\tilde{y}_2)\in \tilde{\Omega}$.
Regarding the boundary conditions, a typical choice is to consider {\itshape strong anchoring } on the boundary of~$\tilde{Q}$. Mathematically, this corresponds to a (non-homogeneous) Dirichlet condition on $\partial\tilde{Q}$, that is $\tilde{U} =G$, where $G$ is a %assigned 
boundary datum. 
% When dealing with these boundary conditions, the 
% mathematical analysis of the Oseen-Frank energy often neglects the last term, the saddle-splay, in the energy density since for this term we have that, as originally observed by C.W. Oseen (see \cite{ericksen76}),
With this choice of the boundary conditions, the saddle-splay term
can often be ignored in the mathematical analysis because, as originally observed 
by C.W. Oseen (see \cite{ericksen76}), there holds
\begin{equation}
\label{eq:NL}
\int_{\tilde{Q}}\left(\textrm{tr}\left(\nabla \tilde{U}\right)^2-(\div\tilde{U})^2\right)\d\tilde{y} 
= \int_{\tilde{Q}}\div \left((\nabla \tilde{U})\tilde{U} -(\div \tilde{U})\tilde{U}\right)\d\tilde{y}.
\end{equation}
By applying the divergence theorem, the right-hand side reduces 
to a surfaces integral on~$\partial\tilde{Q}$ which depends only 
on the Dirichlet datum~$G$ and its tangential derivatives.
In other words, the saddle-splay term is a Null Lagrangian,
which does not affect the minimizers of the Oseen-Frank energy
because it depends {\itshape only} on the (prescribed) 
values of $\tilde{U}$ on $\partial\tilde{Q}$.  

When imposing different type of anchoring (i.e., other boundary conditions),
it may not physically reasonable to neglect the saddle-splay term any longer. 
In our contribution, for example, we impose strong anchoring only on the lateral boundary of the slab, namely on $\partial\tilde{\Omega}\times (0,\tilde{h})$. 
% That is, given $G\colon\partial\tilde{\Omega}\times (0,\tilde{h})\to \mathbb{S}^2$, we prescribe the Dirichlet boundary condition $\tilde{U} = G $ on $\partial\tilde{\Omega}\times (0,\tilde{h})$. 
On the top and on the bottom face of $\tilde{Q}$, we penalize those configurations with non-zero components along the $\hat{e}_3$ direction. We incorporate this boundary condition directly in the energy as a surface energy term, in the form of a Rapini-Papoular {\itshape weak anchoring} condition
(see e.g.~\cite{BarberoDurand}).
Therefore, the energy we are interested in is the following:
\begin{equation}
\label{eq:OF_general}
E(\tilde{U}) = \int_{Q}\sigma(\tilde{U},\nabla\tilde{U}) \, \d \tilde{y} + 
\frac{\tilde{\lambda}}{2}\int_{\Omega\times\left\{0,1\right\}}\abs{(\tilde{U},\nu)}^2 \d \tilde{y}', 
\end{equation}
where $\nu =\pm \hat{e}_3$ denotes the outer normal to~$\tOmega\times \left\{1\right\}$ and to $\tOmega\times \left\{0\right\}$, respectively
and~$\tilde{\lambda}$ is a positive constant that measures the strength of the anchoring. 
% {\BBB We allow that~$\tilde{\lambda}$ might depend on the thickness $h$.
% \textbf{Capisco forse quello che vuoi dire, ma mi pare che 
% questa frase si presti a malintesi (cio\`e', si potrebbe pensare che 
% ci sia proprio una relazione fisica che lega tra loro 
% le due grandezze\ldots).}}

In this paper, we will consider the simplified situation in which the Frank's constant are tuned in such a way that  $k_1=k_2=k_3 = \kappa$ and $k_4=0$. This choice is commonly known as the {\itshape one constant approximation}. 
In this case, the Oseen-Frank energy density reduces to 
\[
2\sigma(\tilde{U},\nabla\tilde{U})= \kappa \abs{\nabla \tilde{U}}^2,
\]
and our energy functional reads 
\begin{equation}
\label{eq:OF_cyl_h}
\tilde{E}(\tilde{U}) = \frac{\kappa}{2}\int_{\tOmega\times (0,\th)}\abs{\nabla \tilde{U}}^2 \d \tilde{y} 
+\frac{\tilde{\lambda}}{2}\int_{\tOmega\times \left\{0,\th\right\}}\abs{(\tilde{U},\nu)}^2 \d \tilde{y}'.
\end{equation}

% {\BBB\textbf{In quello che segue,} faccio una piccola variante: un solo cambio di variabile, anisotropo, anzich\'e due. La versione precedente \`e sotto, commentata.}
Before addressing the mathematical analysis of~\eqref{eq:OF_cyl_h}, 
we write the functional in non-di\-men\-sion\-al form. 
As we plan to study the dimension reduction limit, 
in which the relative thickness of the film tends to zero,
it is convenient to consider an anisotropic rescaling
of the spatial variables, so that the non-dimensional
functional is set on a fixed domain. 
%The director field $z:\tilde{Q}\to \R^3$ belongs to $H^1(\tilde{Q};\R^3)$ with $\abs{z}=1$ almost everywhere. 
%Moreover we ask that $z_{|\partial\tOmega\times (0,h)} = (g(y'), 0)$ with $g\in H^{1/2}(\partial\tOmega,\mathbb{S}^1)$. The degree of $g$ is not necessarily zero. 
%The constant $\tilde{\lambda} \ge 0$ measures the strength of the anchoring condition on the top and on the bottom boundaries $\tOmega\times \left\{1\right\}$ and $\tOmega\times \left\{0\right\}$ ({\color{red}elaborare un po' di più}).
Let~$d_{\tilde{\Omega}}$ be a typical length scale of the two-dimensional cross-section of the domain
--- say, the diameter of~$\tilde{\Omega}$.
We consider rescaled, non-dimensional variables given by
% non dimensionalize the problem by scaling the spatial variables according to 
\[
x_1 = \frac{\tilde{y}_1}{d_{\tOmega}},\qquad x_2 = \frac{\tilde{y}_2}{d_{\tOmega}},\qquad x_3= \frac{\tilde{y}_3}{\tilde{h}}, 
\] 
The variables~$x = (x_1, \, x_2, \, x_3)$ now range in~$Q := \Omega\times (0, \, 1)$, where~$\Omega := \tilde{\Omega}/d_{\tilde{\Omega}}$
is the rescaled, dimension-less domain.
We also introduce the rescaled director field~$U\colon Q\to\R^3$, given by 
$U(x) := \tilde{U}(\tilde{y})$ for $\tilde{y}\in Q$,
and define non-dimensional parameters
\[
 h:= \dfrac{\th}{d_{\tOmega}}, \qquad \lambda:= \dfrac{\tilde{\lambda}d_{\tOmega}}{\kappa}
\]
Then, the Oseen-Frank energy in non-dimensional variables reads as
\begin{equation}
\label{eq:OF_resc1}
F_h(U) := \frac{1}{\kappa \tilde{h}}\tilde{E}(\tilde{U}) = \dfrac{1}{2}\int_{Q}\left(\sum_{k=1}^{2}\abs{\dfrac{\partial U}{\partial x_k}}^2 + \dfrac{1}{h^2}\abs{\dfrac{\partial U}{\partial x_3}}^2\right)\d x + \dfrac{\lambda}{2 h}\int_{\Omega\times\left\{0,1\right\}}\abs{(U,\nu)}^2 \d x'
\end{equation}
for~$x^\prime := (x_1, \, x_2)$.

The asymptotic behavior of $F_h$ depends 
on the interplay between $\lambda$ and $h$. 
In this paper, we assume that
the parameters~$\lambda$ and~$h$ vary along a curve
described by $\lambda=\lambda(h)$, 
and we consider a regime in which
\begin{equation} \label{h,lambda}
 h\to 0, \qquad 
 \frac{\lambda(h)}{h} 
 %= \dfrac{\tilde{\lambda} d_{\tilde{\Omega}}^2}{\kappa\tilde{h}}
 \to +\infty, \qquad
 \sqrt{2 h \, \lambda(h)} \leq 1
\end{equation}
The limit~$h\to 0$ corresponds to dimensional reduction, 
as the thickness of the nematic film
becomes very small compared to the size of two-dimensional cross section.
The conditions on~$\lambda(h)$ quantify the relative 
strength of the anchoring at the boundary,
compared to the other parameters of the problem.

In order to emphasize further the mathematical analogies
between~\eqref{eq:OF_resc1} and the Ginzburg-Landau functional
(as considered by Bethuel, Brezis and H\'elein in~\cite{BBH}),
we define 
\begin{equation}
\label{eq:epsi}
\eps:= \sqrt{\dfrac{h}{\lambda(h)}}.
\end{equation}
Moreover, we assume that the function $h\mapsto\sqrt{\dfrac{h}{\lambda(h)}}$ 
is locally invertible in a (right) neighborhood of $h=0$ and
we define a strictly positive function~$\eta = \eta(\eps)$
as the inverse of $h\mapsto\sqrt{\dfrac{h}{\lambda(h)}}$.
Then, the functional~\eqref{eq:OF_resc1} reduces to~\eqref{eq:energyQ},
and the conditions~\eqref{h,lambda} reduce to~\eqref{eq:ipeta}.

\section{Preliminaries}
\label{sec:prelim}

\subsection{Functional spaces and existence of Oseen-Frank minimizers}
\label{ssec:funct_sp_min}

In this subsection, we recall some classical results that imply
the existence of minimizers for the energy~\eqref{eq:energyQ}. 
From the discussion above, it is clear that the natural 
functional space for the minimization of~\eqref{eq:energyQ}
is the space~$\mathcal{A}_{G}$ defined in~\eqref{eq:funct-frame}.
% where the (lateral) boundary datum $g\in H^{1/2}(\partial\Omega\times (0,1);\S^2)$.
First of all, we observe that this space is nonempty.
We write~$H^{1/2}(\partial\Omega\times (0,1);\S^2)$
for the set of maps~$G\in H^{1/2}(\partial\Omega\times (0, 1); \R^3)$
that satisfy~$\abs{G}=1$ a.e.~on~$\partial\Omega$,
with respect to the surface measure.

\begin{lemma} \label{lemma:nonempty}
For any $G\in H^{1/2}(\partial\Omega\times (0,1);\S^2)$,
we have
\[
 \mathcal{A}_G \neq \emptyset.
\]
\end{lemma}
Lemma~\ref{lemma:nonempty} is actually more general 
than what we need in this context, for it applies 
to boundary data~$G$ that may or may not be of the
form~\eqref{eq:bddatum}. Moreover,
the arguments in this section carry over,
with no substantial change, to prove
existence of minimizers for~$F_\eps$
on general bounded Lipschitz domains $Q\subset \R^3$
such that $\partial Q = \Gamma_1 \cup \Gamma_2$,
with Dirichlet boundary conditions
on~$\Gamma_1$ and weak anchoring conditions on~$\Gamma_2$.
% for it applies 
% to boundary data~$G$
% that may or may not be of the form~\eqref{eq:bddatum} {\BBB and to 
% general bounded Lipschitz domains $Q\subset \R^3$ such that $\partial Q = \Gamma_1 \cup \Gamma_2$. 
% In particular, we might address the minimization of the energy \eqref{eq:energyQ} on this general domain with Dirichlet boundary conditions
% on $\Gamma_1$ and weak anchoring conditions on $\Gamma_2$.  
% }

The proof follows by an argument due to Hardt, Kinderlehrer
and Lin~\cite[Page 556]{HKL86}, which we recall for the reader's convenience.
\begin{proof}[Proof of Lemma~\ref{lemma:nonempty}]
In this proof, contrarily to the rest of the paper, 
we will use the notation~$B_r(0)$ to denote 
\emph{three-dimensional balls}, centered at the
origin, of radius~$r>0$.
% The proof follows the line of \cite[Pag. 556]{HKL86}.
Given a boundary datum $G\in H^{1/2}(\partial\Omega\times(0, \, 1); \, \S^2)$,
we let~$w$ be the harmonic extension of~$G$ to the cylinder~$Q$, namely
\[
w\in \argmin\left\{\frac{1}{2}\int_Q\abs{\nabla v}^2 \d x; \quad v = G\quad \hbox{ on } \partial\Omega\times (0,1)\right\}.
\] 
Elliptic theory implies that~$w$ is smooth inside~$Q$.
Moreover, we have that $\abs{w}\le 1$ in $Q$, by the maximum principle.
Now, if $\min_{Q}\abs{w}>0$, then $V:= \frac{w}{\abs{w}}\in H^1_{G}(Q;\S^2)\neq \emptyset$. 
% In the opposite case, namely $\abs{w}\ge 0$ in $\Omega$,
Otherwise, we show that $H^1_{G}(Q;\S^2)\neq \emptyset$
by reasoning as follows. 

First of all, we let $z\in B_{1}(0)\subset\R^3$ be a regular value of $w$
(that is, a value~$z$ such that~$\det\nabla w(x)\neq 0$ 
for any~$x\in B_1(0)$ with~$w(x) = z$).
By Sard's Theorem, we know that the set of regular values is dense in $\R^3$.
Moreover, %since $w$ is smooth, we have that  
the inverse function theorem implies that
$Z= \left\{x\in Q: w(x) = z\right\}$ is a discrete set
(with respect to the Euclidean topology on~$Q$),
hence locally finite in $Q$. 

We consider the map $v(x,z) := \frac{w(x)-z}{\abs{w(x)-z}}$,
which is smooth in the variable~$x\in Q\setminus Z$ and takes values in~$\S^2$.
% We claim that
% \begin{equation} \label{nonempty1}
%  v(\cdot, \, z)\in H^1_{\mathrm{loc}}(Q; \, \S^2)
% \end{equation}
% Let~$A\subset Q$ be an open set with $\bar{A}$ compact in~$Q$. 
% As~$Z$ is locally finite, we can write $Z\cap A= \bigcup_{k=1}^N \left\{x_k\right\}$ for some finite~$N$. The map $v$ is clearly smooth in $Q\setminus \bigcup_{k=1}^N B_{\delta}(\bar{x}_{k})$, where $\delta > 0$ is chosen so small in such a way that the balls $B_\delta(\bar{x}_k)$ are disjoint. Therefore, 
% in order to prove~\eqref{nonempty1}, it suffices to check that $v(\cdot, \, z)\in H^1(B_{\delta}(\bar{x}_k); \, \S^2)$ for any $k=1,\ldots,N$.
% Fix $k\in \left\{1,\ldots,N\right\}$ and expand $w$ around $\bar{x}_k$.
% We have 
% \[
% w(x) - z= \nabla w(\bar{x}_k)(x-\bar{x}_k) + \mathrm{o}(\abs{x-\bar{x}_k})
% \qquad \textrm{as } x\to x_k. 
% \]
% and hence (recall that~$\nabla w(\bar{x}_k)$ is non-singular, because~$z$
% is a regular value),
% \[
% v(x, \, z) = \frac{\nabla w(\bar{x}_k)(x-\bar{x}_k)}{\abs{\nabla w(\bar{x}_k)(x-\bar{x}_k)}} + \mathrm{o}(\abs{x-\bar{x}_k})
% \qquad \textrm{as } x\to x_k. 
% \]
% A simple computation show that
% $\frac{M(x-\bar{x}_k)}{\abs{M(x-\bar{x}_k)}}\in H^1(B_\delta(\bar{x}_k);\S^2)$ for any non singular $3\times 3$ matrix $M$. 
% Consequently, we have proved~\eqref{nonempty1}.
% % obtained that $v\in H^1_{\textrm{loc}}(Q;\S^2)$. 

We claim that, upon choosing a suitable~$z$ among the regular values of $w$, there holds
\begin{equation} \label{nonempty2}
\nabla v(\cdot, \, z) \in L^2(Q;\R^{3\times 3}).
\end{equation}
Here~$\nabla v$ denotes the classical gradient of~$v$
with respect to the variable~$x$, which is well-defined
everywhere in~$Q\setminus Z$.
Indeed, for any $x\in Q\setminus Z$, we have that 
\[
\nabla v(x,z) = \frac{1}{\abs{w(x)-z}}\left(\nabla w(x) -\frac{1}{\abs{w(x)-z}^2}(w(x)-z) \otimes (\nabla w(x))(w(x)-z)\right) \! ,
\]
and hence,
\begin{equation} \label{nonempty2.5}
\abs{\nabla v(x, \, z)}\le \frac{\abs{\nabla w(x)}}{\abs{w(x)-z}}.
\end{equation}
We show that 
\begin{equation} \label{nonempty3}
\int_{B_1(0)}\abs{w(x)-z}^{-2}\d z <+\infty \quad \textrm{ for almost any }x\in Q.
\end{equation}
Indeed, setting $y:= w(x)-z$ and recalling that~$\abs{w(x)} \leq 1$ for almost e any~$x\in Q$, we readily obtain 
\[
\int_{B_1(0)}\abs{w(x)-z}^{-2}\d z \leq \int_{B_2(0)}\abs{y}^{-2}\d y = 8\pi \quad \textrm{ for almost any }x\in Q,
\]
which proves~\eqref{nonempty3}. As a consequence, we obtain
\begin{equation} \label{nonempty4}
\int_{Q}\abs{\nabla w(x)}^2\left(\int_{B_1(0)}\abs{w(x)-z}^{-2}\d z\right)\d x \leq 8\pi\int_{Q}\abs{\nabla w(x)}^2\d x<+\infty.
\end{equation}
Therefore, thanks to~\eqref{nonempty2.5}, \eqref{nonempty4} 
and Tonelli's Theorem (see, e.g., \cite[Theorem 4.4]{brezis}), 
we conclude that 
\[
(z,x)\mapsto \nabla v(x,z)\in L^2(Q\times B_1(0);\R^{3\times 3})
\]
and that
\[
 \begin{split}
%   \int_{Q\times B_1(0)}\abs{\nabla v(x,z)}^2\d x\d z = \int_{Q}\left(\int_{B_1(0)}\abs{\nabla v(x,z)}^2 \d z\right)\d x 
%   = \int_{B_1}\left(\int_{Q}\abs{\nabla v(x,z)}^2 \d x\right)\d z <+\infty,
 \int_{B_1}\left(\int_{Q}\abs{\nabla v(x,z)}^2 \d x\right)\d z
  = \int_{Q}\left(\int_{B_1(0)}\abs{\nabla v(x,z)}^2 \d z\right)\d x 
  \leq 8\pi\int_{Q}\abs{\nabla w(x)}^2\d x,
 \end{split}
\]
which implies that we can choose $z$ among the regular values of $w$ in such a way that 
\[
\int_{Q}\abs{\nabla v(x,z)}^2 \d x
\leq 6\int_{Q}\abs{\nabla w(x)}^2\d x<+\infty, 
\]
thus proving~\eqref{nonempty2}.

We must check that~$\nabla v(\cdot, z)$ coincides with 
the gradient in the sense of distributions on~$Q$.
To this end, take a test
function~$\varphi\in C^\infty_{\mathrm{c}}(Q, \, \R^3)$.
The set~$Z$ is locally finite, so only finitely many points of~$Z$
are contained in the support of~$\varphi$. Then, using 
cut-off functions, we can construct a
sequence of functions~$\psi_n\in C^\infty(Q, \, \R)$ 
such that $\psi_n = 0$ in a neighbourhood of each point
of~$Z\cap \mathrm{support}(\varphi)$ and~$\psi_n\to 1$ 
strongly in~$H^1(Q)$. As~$v(\cdot, z)$ is smooth on the support 
of~$\psi_n \varphi$, integrating by parts gives
\[
 \int_Q \psi_n(x) \partial_k v(x, z)\cdot \varphi(x) \, \d x
 = -\int_Q \partial_k\psi_n(x) \, v(x, z)\cdot \varphi(x) \, \d x
  -\int_Q \psi_n(x) v(x, z)\cdot \partial_k\varphi(x) \, \d x
\]
for~$k\in\{1, \, 2, \, 3\}$.
By passing to the limit as~$n\to+\infty$, we conclude that
that~$\nabla v(\cdot, z)$ is indeed the distributional gradient
of~$v$ and hence, thanks to~\eqref{nonempty2},
that~$v(\cdot, z)\in H^1(Q; \S^2)$.

To conclude, we have to adjust the boundary conditions.
% This is done as in \cite{HKL86}.

Indeed, the map~$v(\cdot, z)$ does not agree with~$G$
on~$\partial\Omega\times (0, 1)$; instead, it satisfies
\begin{equation} \label{nonempty5}
 v(x, z) = \frac{G(x) - z}{\abs{G(x) - z}}
 \qquad \textrm{for a.e. }
 x\in \partial\Omega\times (0, \, 1),
\end{equation}
in the sense of traces. Now, let us consider the
map~$\Pi_z\colon\S^2\to\S^2$ given by
\[
 \Pi_z(y) := \frac{y - z}{\abs{y - z}}
 \qquad \textrm{for any } y\in \S^2.
\]
This map is well-defined and smooth, for we have 
chosen~$z$ such that~$\abs{z} < 1$. 
An explicit computation shows that the differential of~$\Pi_z$
is non-singular at each point. Therefore, the local inversion
theorem implies that~$\Pi_z$ is a local diffeomorphism.
On the other hand, it is not hard to check that~$\Pi_z$ is one-to-one;
then, since~$\S^2$ is compact and connected,
$\Pi_z$ must be a global diffeomorphism $\S^2\to\S^2$.
Therefore, the
$V := \Pi_z^{-1}\circ v(\cdot, z)$ belongs to~$H^1(\Omega; \S^2)$
and, thanks to~\eqref{nonempty5}, it satisfies~$V = G$ 
on~$\partial\Omega\times (0, \, 1)$, in the sense of traces.
\end{proof}

Once we have proved that $\mathcal{A}_{G}\neq \emptyset$, the next result follows from a routine application of the Direct Method of the Calculus of Variations. 
\begin{prop} \label{prop:exists}
For any $G\in H^{1/2}(\partial \Omega\times (0,1);\S^2)$ and any $\eps>0$, there exists 
\[
U_\eps\in \Argmin\left\{ F_\eps(U): U\in \mathcal{A}_G\right\}.
\]
\end{prop}

\subsection{Distributional Jacobian}
\label{ssec:jaco}
We recall the definition of the Jacobian determinant in the sense of distributions, following~\cite[Chapters 1 \& 2]{BM21}.
For any $u\in W^{1,1}_{\mathrm{loc}}(\Omega; \R^2)\cap L^{\infty}_{\mathrm{loc}}(\Omega;\R^2)$, we define the vector field
$j(u)\colon\R^2\to\R^2$ as 
\begin{equation}
\label{eq:prejaco}
j(u) := u_1 D u_2 - u_2 D u_1,
\end{equation}
which in components reads 
\[
(j(u))_k = u_1\frac{\partial u_2}{\partial x_k} - u_2 \frac{\partial u_1}{\partial x_k} = u\times \frac{\partial u}{\partial x_k}
\qquad \textrm{for }k=1,2.
\]
We observe that $j(u)\in L^1_{\mathrm{loc}}(\Omega;\R^2)$, by construction. Therefore, it makes sense to define the distribution 
\begin{equation}
\label{eq:jaco}
\langle J u,\phi\rangle : = -\int_{\Omega} (j(u), \nabla^{\perp} \phi) \, \d x' \qquad \textrm{for any } \phi\in C^{0,1}_c(\Omega)
\end{equation}
(where 
$\nabla^{\perp}\phi:= \left(-\frac{\partial \phi}{\partial x_2}, \frac{\partial \phi}{\partial x_1}\right)$), namely 
\[
J u = \textrm{curl} \, j(u)
\]
in the sense of distributions. 
We observe that if~$u$ is smooth, then 
\begin{equation}
\label{eq:smooth_jaco}
Ju = 2 \, \textrm{det}D u.
\end{equation}
The (distributional) Jacobian can be related to the
notion of degree for maps $u\colon U\to \R^2$,
where $U$ is a bounded, regular (say, Lipschitz) domain 
$U\subset \R^2$. Suppose that~$u$ is sufficiently smooth on~$\partial U$
and satisfies~$\abs{u}=1$ on $\partial U$.
We define 
%and suppose that $u\in W^{1,1}(U,\R^2)\cap L^{\infty}(U;\R^2)$, then we define (the integer)
\begin{equation}
\label{eq:degree_def}
\textrm{deg}\left(u,\partial U\right):= 
\frac{1}{2\pi}\int_{\partial U}j(u)\cdot \tau \, \d s
= \frac{1}{2\pi} \int_{\partial U} u\times \partial_\tau u \, \d s, 
\end{equation} 
where $\tau$ is the unit tangent to~$\partial U$, positively 
oriented with the outer normal, and~$\partial_\tau u$ 
is the tangential derivative of~$u$. 
Note that, thanks to \eqref{eq:smooth_jaco}, 
we have (formally at least)
\begin{equation}
\label{eq:degree_Jac}
\textrm{deg}(u, \partial U) = \frac{1}{\pi}\int_{U}\det D u \, \d x.
\end{equation}
Moreover, the degree is an integer number. 
Indeed, suppose (for simplicity only) that~$\partial U$
is connected, i.e.~$U$ is simply connected.
By parametrizing~$\partial U$ by a Lipschitz curve 
$\gamma\colon [0, \, 1]\to\partial U$ and writing~$u(\gamma(t)) = \exp(i\theta(t))$ for some Lipschitz function~$\theta\colon [0, \, 1]\to\R$, 
we obtain
\[
 \textrm{deg}\left(u,\partial U\right)
  \stackrel{\eqref{eq:degree_def}}{=}
  \frac{1}{2\pi} \int_0^{2\pi} \theta^\prime(t) \, \d t
  = \frac{1}{2\pi} \left(\theta(2\pi) - \theta(0)\right) 
\]
which proves that~$\deg(u,\partial U)$ is an integer,
because $\exp(i\theta(0))= u(\gamma(0)) = u(\gamma(1)) = \exp(i\theta(1))$.

Equations~\eqref{eq:degree_def} and~\eqref{eq:degree_Jac}
make sense when~$u$ is regular enough. 
Boutet de Monvel and Gabber~\cite[Appendix]{BoutetdeMonvel}
observed that~$\deg(u,\partial U)$ is well-defined
as soon as~$u\in H^{1/2}(\partial U; \S^1)$,
because the right-hand side of~\eqref{eq:degree_def} 
can then be interpreted as a duality pairing between~$u\in H^{1/2}(\partial U)$
and~$\partial_\tau u\in H^{-1/2}(\partial U)$.
Equation~\eqref{eq:degree_Jac} holds true, for instance,
when~$u\in H^1(U;\R^2)$, for then the Jacobian~$J u$
is an integrable function and the identity~\eqref{eq:smooth_jaco}
is satisfied pointwise almost everywhere, 
as shown by a density argument. However,
Equation~\eqref{eq:degree_Jac} remains valid
in other cases, too --- for instance, when~$u\in H^{1/2}(\partial U;\S^1)$
and~$J u$ is a measure, which may be the case even if~$u\notin H^1(U;\R^2)$.

It is interesting to remark that \eqref{eq:degree_def} is one of the possible (albeit equivalent) definitions of degree. 
In particular, the definition in terms of the Jacobian determinant is particularly suited for the applications to the Ginzburg-Landau energy. 
For other definitions we refer to \cite{BM21}.

\subsection{Canonical harmonic map}
\label{ssec:chm}
We let $g:\partial\Omega\to \mathbb{S}^1$ be an assigned boundary datum with 
\[
g\in H^{1/2}(\partial\Omega;\mathbb{S}^1).
\]
We recall a definition from~\cite{BBH}.
\begin{definition}[Canonical Harmonic map]
\label{def:chm}
Given $N$ distinct points $\a=(a_1,\ldots,a_N)$ in $\Omega$ and integers $\db=(d_1,\ldots,d_N)$ such that  
\begin{equation}
\label{eq:PH}
\sum_{n=1}^N d_n = \textrm{deg}(g,\partial\Omega), 
\end{equation}
we say that $u^*\colon\Omega\to\S^1$ is a {\itshape Canonical Harmonic Map}
associated with $(\a, \, \db, \, g)$ if $u^*\in W^{1,p}(\Omega;\mathbb{S}^1)\cap H^{1}_{\mathrm{loc}}\left(\Omega\setminus \bigcup_{n=1}^N\left\{a_n\right\};\mathbb{S}^1\right)$ for any $p\in [1,2)$, if $u^* = g$
on $\partial \Omega$ and if it satisfies
\begin{equation}
\begin{split}
\label{eq:chm}
&\div j(u^*) = 0 \qquad \textrm{ in } \Omega,\\
&Ju^* = 2\pi \sum_{n=1}^N d_n\delta_{a_n}.
\end{split}
\end{equation}
\end{definition}

If~$d_j = 1$ for any index~$j$
(as it happens in Theorem~\ref{th:minimizers}, for instance), 
we omit the dependence on~$\db$ and say that~$u^*$
is a Canonical Harmonic Map associated with~$\a$ and~$g$.
Once we fix the boundary datum and the collection of points and integers $(\a,\db)$ satisfying \eqref{eq:PH}, there exists
a unique Canonical Harmonic map $u^*$ for $(\a,\db,g)$ (see \cite[Section I.3]{BBH}). 
Being $j(u^*)$ divergence-free in a simply connected domain $\Omega$, the Hodge decomposition theorem entails that there exists $\Psi\colon\Omega\to \R$ such that 
\begin{equation}
\label{eq:hodge}
j(u^*) = \nabla^{\perp}\Psi,
\end{equation}
namely
\[
\begin{cases}
\displaystyle u^*\times \frac{\partial u^*}{\partial x_1} = -\frac{\partial \Psi}{\partial x_2}\qquad &\textrm{in }\Omega\\[.4cm]
\displaystyle u^* \times \frac{\partial u^*}{\partial x_2} = \frac{\partial \Psi}{\partial x_1}\qquad &\textrm{in }\Omega.
\end{cases}
\]
Therefore, as a distribution, we have
\[
J u^*= \textrm{curl} \nabla^{\perp}\Psi = \Delta \Psi
\]
and thus the second equation in~\eqref{eq:chm} and the boundary condition~$u^* = g $ on $\partial \Omega$ imply that 
\begin{equation}
\begin{cases}
\label{eq:poisson}
\displaystyle\Delta \Psi = 2\pi \sum_{n=1}^N d_n\delta_{a_n} \qquad &\textrm{in }\Omega\\[.4cm]
\displaystyle\frac{\partial\Psi}{\partial \nu_\Omega} = g\times \partial_\tau g \qquad &\textrm{on }\partial\Omega,
\end{cases}
\end{equation}
where $\nu_\Omega$ is the outer normal to $\partial\Omega$ and $\tau$ is the tangent unit vector to~$\partial\Omega$, positively oriented with respect to $\nu_\Omega$. 

\subsection{Renormalized Energy}
\label{sse:reno}
Given a boundary datum $g\in H^{1/2}(\partial\Omega;\mathbb{S}^1)$,  $N$ distinct points $\a =(a_1,\ldots,a_N)$ in $\Omega$ and integers $\db =(d_1,\ldots,d_N)$ that satisfy~\eqref{eq:PH},
% such that  
% \[
% \sum_{n=1}^N d_n = \textrm{deg}(g,\partial\Omega), 
% \]
let~$u^*$ be the unique Canonical Harmonic map associated to~$(\a,\db,g)$. 
We define the Renormalized Energy of~$(\a, \db)$ as 
\begin{equation}
\label{eq:RE_def}
W_g(\a, \db):= \lim_{\sigma\to 0}\left(\frac{1}{2}\int_{\Omega\setminus \bigcup_{n=1}^N B_{\sigma}(a_n)}\abs{ D u^*}^2 \d x' + \pi \left(\sum_{n=1}^N d_n^2\right)\log \sigma\right).
\end{equation}
The above limit exists according to \cite[Theorem I.8]{BBH}. 
Moreover, the Renormalized Energy is a smooth function 
%$W_g:\Omega^N\times \mathbb{Z}^N\to \R$ 
can be equivalently written (see \cite[Theorem I.7]{BBH}) as
\begin{equation}
\label{eq:RE}
W_g(\a, \db) = -\pi \sum_{i\neq j}^{N}d_i d_j \log\abs{a_i-a_j} + 
\frac{1}{2}\int_{\partial \Omega}\Psi\left(g\times g_\tau\right) -\pi \sum_{i=1}^{N} d_i R (a_i),
\end{equation}
where~$\Psi$ is the unique solution of~\eqref{eq:poisson}
that satisfies
\[
 \int_{\partial\Omega} \Psi \,\d s = 0
\]
and~$R\colon\Omega\to \R$ is the smooth (harmonic) function defined by
\[
 R(x) := \Psi (x) - \sum_{j=1}^N d_j \log\abs{a-a_j}
 \qquad \textrm{for } x\in \Omega\setminus\{a_1, \ldots, a_N\}
\]
and extended to~$\Omega$ by continuity.

When~$d_i=1$ for any $i=1,\ldots,N$ (as it happens in Theorem~\ref{th:minimizers}, for instance), we omit 
the dependence on~$\d$ and write~$W_g(\a)$ 
instead of~$W_g(\a, \, \db)$. We have 
\begin{equation}
\label{eq:boundary_coale}
W_g(\a)\rightarrow +\infty \qquad \textrm{ if }\qquad \min\left\{\min_{i\neq j}\abs{a_i-a_j}, \min_i \textrm{dist}\left(a_i,\partial \Omega\right)\right\}\to 0.
\end{equation}

%%%%%%%%%%%%%%%%%%%%%%%

\section{The $\Gamma$-convergence result}
\label{sec:proofGamma}

\subsection{A preliminary computation: emergence of a Ginzburg-Landau energy}
\label{ssec:prelim_compu}

Let~$U\in \mathcal{A}_{G}$. 
% First of all we notice that 
% \begin{equation}
% \label{eq:bound1}
% \frac{1}{2}\int_{Q}\abs{\frac{\partial U}{\partial x_3}}^2\d x \le \eta^2(\eps)F_\eps(U).
% \end{equation}
We decompose~$U$ as in \eqref{eq:ortogonal_decoU}, namely
\begin{equation}
\label{eq:decoU}
U = \Pi(U) + \left(U,\hat{e}_3\right)\hat{e}_3 = u + U^{\perp}\hat{e}_3,
\end{equation}
and we let 
\begin{equation}
\label{eq:mean}
\bar{u}:= \int_{0}^1 \Pi(U)\, \d x_3 = \int_{0}^1 u \, \d x_3.
\end{equation}
%%%
We prove the following:

\begin{lemma}
\label{eq:lemub}
Let $U\in \mathcal{A}_{G}$. 
Then
\[
\norm{\bar u-u(\cdot, x_3)}_{L^2(\Omega)} 
\leq \frac{1}{\sqrt{2}} \norm{\frac{\partial U}{\partial x_3}}_{L^2(Q)}
% \le \eta(\eps)\sqrt{F_\eps(U)} 
\qquad \textrm{ for a.e. } x_3\in (0,1).
\]
\end{lemma}
\begin{proof}
% We prove the Lemma with $x_3 =0$ just for the sake of notational simplicity. 
First of all we observe that for almost any $x'\in \Omega$
and~$x_3\in (0, 1)$, there holds
\[
\begin{split}
\bar{u}(x')-u(x',x_3) = \int_{0}^1 \left(u(x',t) - u(x',x_3)\right)\d t = \int_{0}^1 \left(\int_{x_3}^t\frac{\partial u(x',s)}{\partial x_3}\d s\right) \d t.
\end{split}
\]
Thus, using Jensen's inequality twice and integrating by parts,
% (separately on the intervals~$(0, \, x_3)$ and~$(x_3, \, 1)$), 
we obtain
\[
\begin{split}
\abs{\bar{u}(x')-u(x',x_3)}^2 &= 
\left(\int_{0}^1 \left(\int_{x_3}^t\frac{\partial u(x',s)}{\partial x_3}\d s\right) \d t\right)^2 
\le \int_{0}^1 \abs{t - x_3}\abs{\int_{x_3}^t \abs{\frac{\partial u(x',s)}{\partial x_3}}^2\d s} \d t\\
% &  = \frac{1}{2}\int_{0}^1 (1-t^2) \abs{\frac{\partial u(x',t)}{\partial x_3}}^2 \d t.
&  = \frac{1}{2}x_3^2 \int_0^{x_3} \abs{\frac{\partial u(x',t)}{\partial x_3}}^2 \d t
+ \frac{1}{2}(1 -x_3)^2 \int_{x_3}^1 \abs{\frac{\partial u(x',t)}{\partial x_3}}^2 \d t\\
&\hspace{2cm} - \frac{1}{2}\int_{0}^1 (x_3-t)^2 \abs{\frac{\partial u(x',t)}{\partial x_3}}^2 \d t \\
&\leq \frac{1}{2}\int_0^1 \left(x_3^2 + (1 - x_3)^2 - (x_3 - t)^2\right) \abs{\frac{\partial u(x',t)}{\partial x_3}}^2 \d t
\end{split}
\]
A routine computation 
shows that $0 \leq x_3^2 + (1 - x_3)^2 - (x_3 - t)^2 \leq 1$
for any~$(x_3, t)\in [0, 1]^2$.
Consequently, integrating with respect to~$x^\prime \in\Omega$,
the lemma follows.
% we obtain
% \[
% \int_{\Omega}\abs{\bar{u}(x')-u(x',0)}^2 \d x' 
% \le \frac{1}{2}\int_{Q} \abs{\frac{\partial u(x',t)}{\partial x_3}}^2 \d x' \d t,
% \]
% which readily implies the thesis, in view of~\eqref{eq:bound1}.
% % recalling that~\eqref{eq:bound1} implies
% % \[
% % \frac{1}{2}\int_{Q}\abs{\frac{\partial u(x',t)}{\partial x_3}}^2 \d x' \d t\le \eta^2(\eps)F_\eps(U).
% % \]
\end{proof}

We rewrite the surface integral in the energy~$F_\eps$ 
as a Ginzburg-Landau type penalization. 
Recalling the orthogonal decomposition~\eqref{eq:decoU}
and that $\abs{U}^2=1$, we have (as $\nu = \pm \hat{e}_3$ on $\Omega\times\left\{1\right\}$ and on $\Omega\times\left\{0\right\}$, respectively)
\begin{equation} 
\label{eq:lb0}
\frac{1}{\eps^2}\int_{\Omega\times\left\{0,1\right\}}\abs{(\nu,U)}^2 \d x' = 
\frac{1}{\eps^2}\int_{\Omega\times\left\{0,1\right\}} \left(1-\abs{u}^2\right)\d x',
\end{equation}
Note that the integrand at the right-hand side 
is non negative, since $\abs{u}^2\le 1$.
Since 
\[
0\le (1-x^2)^2\le 1-x^2 \quad \textrm{ for } \abs{x}\le 1,
\]
we have 
\begin{equation}
\label{eq:lb1}
\int_{\Omega}\left(1-\abs{u(x',0)}^2\right)\d x' \ge \int_{\Omega}\left(1-\abs{u(x',0)}^2\right)^2 \d x'.
\end{equation}
Moreover %(recall that $\abs{\bar u}\le 1$ by construction), 
\begin{equation}
\label{eq:lb2}
\int_{\Omega}\left(1-\abs{\bar{u}}^2\right)^2 \d x' \le 2\left(\int_{\Omega}\abs{\abs{u(x',0)}^2-\abs{\bar{u}(x')}^2}^2 \d x'  +  \int_{\Omega}\left(1-\abs{u(x',0)}^2\right)^2\d x' \right).
\end{equation}
We concentrate on the first term in the right hand side. 
Recalling that $\abs{u}\leq 1$, $\abs{\bar{u}}\leq 1$
by construction and using Lemma \ref{eq:lemub}, we obtain
\[
\begin{split}
\int_{\Omega}\abs{\abs{u(x',0)}^2-\abs{\bar{u}(x')}^2}^2 \d x' &= \int_{\Omega}\abs{\left(u(x',0)-\bar{u}(x'),u(x',0)+\bar{u}(x') \right)}^2\d x' \\
&\le 4 \norm{\bar{u}(\cdot)-u(\cdot,0)}^2_{L^2(\Omega)} 
\le 2 \int_Q \abs{\frac{\partial U}{\partial x_3}}^2 \d x
\end{split}
\]
(Note that we can take~$x_3 = 0$ in the statement of Lemma~\ref{eq:lemub},
so long as we interpret~$u(\cdot, 0)$ in the sense of traces.)
Therefore, taking~\eqref{eq:lb0}, \eqref{eq:lb1} 
and~\eqref{eq:lb2} into account, we deduce
\begin{equation} 
\label{eq:GLbound-potential}
 \frac{1}{4\eps^2}\int_{\Omega}\left(1-\abs{\bar{u}}^2\right)^2 \d x'
 \leq \frac{1}{\eps^2} \int_Q \abs{\frac{\partial U}{\partial x_3}}^2 \d x
 + \frac{1}{2\eps^2}\int_{\Omega\times\left\{0,1\right\}}\abs{(\nu,U)}^2 \d x'
%  \leq \left(1 + \frac{2\eta^2(\eps)}{\eps^2} \right)  F_\eps(U) 
\end{equation}
%We have
%\[
%\int_{\Omega}\left(1-\abs{\bar{u}}^2\right) \d x'  = 
%\int_{\Omega}\left(\abs{\bar{u}(x')}^2-\abs{u(x',0)}^2\right) \d x' + \int_{\Omega}\left(1-\abs{u(x',0)}^2\right)\d x',
%\]
%and (recalling that $\abs{\Omega}\le 1$ and Lemma \ref{eq:lemub})
%\[
%\begin{split}
%\abs{\int_{\Omega}\left(\abs{\bar{u}(x')}^2-\abs{u(x',0)}^2\right) \d x' }&\le 
%\int_{\Omega}\abs{\abs{\bar{u}(x')}^2-\abs{u(x',0)}^2} \d x'\\
% &= \int_{\Omega}\abs{\left(\bar{u}(x')-u(x',0),\bar{u}(x')+u(x',0)\right)}\d x' \\
%&\le 2\norm{\bar{u}(\cdot)-u(\cdot,0)}_{L^2(\Omega)}\le 2\eta(\eps)\sqrt{F_\eps(U)}.
%\end{split}
%\]
On the other hand, %since $\abs{\nabla_\eps U}^2 \ge \abs{D u}^2$,
Jensen's inequality and Fubini theorem give that 
\begin{equation}
\label{eq:Jenses}
% F_\eps(U)\ge 
\frac{1}{2}\int_{Q}\abs{D u}^2\d x = \frac{1}{2}\int_{\Omega}\left(\int_{0}^1 \abs{D u}^2\d x_3\right)\d x'\ge \frac{1}{2}\int_{\Omega}\abs{D \bar u}^2 \d x'.
\end{equation}
As a result, we can control the Ginzburg-Landau
energy of~$\bar{u}$ in terms of~$F_\eps(U)$, namely 
\begin{equation}
\label{eq:GLbound1}
GL_\eps(\bar{u}):= \frac{1}{2}\int_{\Omega}\abs{D \bar u}^2 \d x' + \frac{1}{4\eps^2}\int_{\Omega}\left(1-\abs{\bar{u}}^2\right)^2 \d x'
\le \max\left(1, \, \frac{2\eta^2(\eps)}{\eps^2}\right)
 F_\eps(U).
\end{equation}
In particular, in the regime of parameters~\eqref{eq:ipeta}
we are interested in,
% --- that is, $\eta(\eps)\ll\eps \ll 1$ ---
we have
\begin{equation}
\label{eq:GLbound2}
GL_\eps(\bar{u}) \leq F_\eps(U)
\end{equation}
for any~$U\in\mathcal{A}_G$.

\begin{remark} \label{rk:morepreciseGLbound}
 %{\BBB Non so se servir\`a! Si pu\`o togliere. }
 If we replace the assumption~\eqref{eq:ipeta}
 with a slightly stronger condition, we 
 can obtain a slightly better bound. Namely,
 assume that there exists a number~$c_* > 0$ 
 such that $2\eta(\eps)^2 \leq (1 - c_*) \eps^2$
 for any~$\eps>0$. Then, for any~$U\in\mathcal{A}_G$
 there holds
 \[
  GL_\eps(\bar{u}) + \frac{c_*}{2\eta(\eps)^2} 
  \int_Q \abs{\frac{\partial U}{\partial x_3}}^2 \d x
  \leq F_\eps(U) 
 \]
 This inequality follows from~\eqref{eq:GLbound-potential}
 and~\eqref{eq:Jenses}.
\end{remark}

% {\color{red}Per avere davvero il $4$ forse serve mettere un'altra costante di fronte al pezzo di penalizzazione nell'energia $F_\eps$ (vedere \eqref{eq:energyQ})}
% 
% {\color{blue}
% Following the discussion at the end of the paragraph on scaling. If we work with $\lambda_h$ satisfying \eqref{eq:choicelambda} then we obtain the following informations on the second term in the right hand side. 
% If $\alpha\in [0,1)$ then we are in the ''good'' regime, namely $\eta(\eps)$ satisfies \eqref{eq:ipeta}, and we obtain that if $U\in \mathcal{A}_{G}$ is chosen in such way that $F_\eps(U)\le C\abs{\log\eps}$ for some constant $C>0$
% \[
% \dfrac{\eta^2(\eps)}{\eps^2} F_\eps(U) = o_{\eps\to 0}(1),
% \]
% }
% and thus
% we have that 
% \begin{equation}
% \label{eq:GLbound}
% \frac{1}{2}\int_{\Omega}\abs{D \bar u}^2 \d x' + \frac{1}{4\eps^2}\int_{\Omega}\left(1-\abs{\bar{u}}^2\right)^{{\color{blue}2}} \d x'\le  F_\eps(U) +  o_{\eps\to 0}(1).
% \end{equation}
% {\color{blue}Mi sembra che ora questa parte qui sotto non serve più. Se non erro, possiamo permetterci di prendere anche $\alpha=0$ (Olè olè olè!!!!)
% 
% On the contrary, if $\alpha=0$ (namely $h = \eps^2$), then if $U\in \mathcal{A}_{G}$ is chosen in such way that $F_\eps(U)\le C_1\abs{\log\eps}$ for some constant $C_1>0$,
% then
% \begin{equation}
% \label{eq:GLbound2}
% \frac{1}{2}\int_{\Omega}\abs{D \bar u}^2 \d x' + \frac{1}{4\eps^2}\int_{\Omega}\left(1-\abs{\bar{u}}^2\right) \d x'\le  C_2 F_\eps(U),
% \end{equation}
% with $C_2$ possibly larger that $C_1$. }

%%%%%%%%%%%%%%%%%%%%%%%%

\subsection{Core Energy}
\label{ssec:core}

For any~$\sigma > 0$, we set $C_{\sigma}:=B_\sigma\times (0,1)$
(where, we recall, $B_\sigma\subset\R^2$
is a two-dimensional disk centered at the origin)
and we consider the minimization problem 
\begin{equation}
\label{eq:core}
\begin{split}
\gamma(\sigma,\eps) &:= \min\bigg\{\frac{1}{2}\int_{C_\sigma}\abs{\nabla_\eps V}^2 \d x + \frac{1}{2\eps^2}\int_{B_\sigma\times \left\{0,1\right\}}\abs{(V,\nu)}^2\d x', \\
&\displaystyle \hspace{3cm} V\colon C_{\sigma}\to \mathbb{S}^2, \quad V = H \quad \textrm{ on }\partial B_{\sigma}\times (0,1) \bigg\},
\end{split}
\end{equation}
where $H(x) := \left(\frac{x'}{\abs{x'}},0\right)$ for any 
$x = (x', x_3)\in \left(B_\sigma\setminus \left\{0\right\}\right)\times (0,1)$. 
The existence of minimizers for~\eqref{eq:core}
follows from Proposition~\ref{prop:exists}.
We have the following (cf.~\cite[Lemma III.1]{BBH})

\begin{lemma}[Core energy]
\label{lemma:core}
The function
\begin{equation*}
 (\sigma, \eps) \mapsto \gamma(\sigma,\eps) - \pi\log\frac{\sigma}{\eps}
\end{equation*}
is bounded and, for any given value of~$\eps$,
is monotonically nonincreasing as a function of~$\sigma$.
Moreover, for any sequence~$\eps_k\to 0$ there exists
a (non-relabelled) subsequence such that the limit
\begin{equation}
\label{eq:conv_core}
\begin{split}
\gamma := \lim_{\sigma\to 0} \, \lim_{k\to+\infty}
 \left(\gamma(\sigma,\eps_k) - \pi\log\frac{\sigma}{\eps_k}\right)
\end{split}
\end{equation}
exists and is finite.
\end{lemma}
In the context of the Ginzburg-Landau functional,
the number~\eqref{eq:conv_core} is sometimes known as the core energy.
However, while the value of~$\gamma$ is uniquely determined
in the Ginzburg-Landau case, here we do not have investigated whether the value of~$\gamma$ 
might depend on the subsequence~$\eps_k\to 0$.

\begin{proof}[Proof of Lemma~\ref{lemma:core}]
For simplicity of notation, let
\[
 \tilde{\gamma}(\sigma, \, \eps)
 := \gamma(\sigma, \eps) - \pi\log\frac{\sigma}{\eps}
\]
We split the proof into steps.

\setcounter{step}{0}
\begin{step}[$\tilde{\gamma}$ is decreasing with respect to~$\sigma$]
 
 Let~$\eps > 0$, $0 < \sigma_1 < \sigma_2$ be given.
 Let~$U_1\in H^{1}(C_{\sigma_1};\mathbb{S}^2)$ be such that 
 \[
  \gamma(\sigma_1, \eps) = F_{\eps}\left(U_1,C_{\sigma_1}\right).
 \]
 Then we consider the field 
 \[
 U:= \begin{cases}
  U_1 &\qquad \textrm{ in } C_{\sigma_1}\\
  H &\qquad \textrm{ in } C_{\sigma_2}\setminus C_{\sigma_1}.
 \end{cases}
 \]
 The field $U$ is, by construction, a competitor for the minimization problem related to~$\gamma(\sigma_1, \eps)$. 
 Therefore, observing that $H\equiv \Pi(H)$ 
 in $\left(B_{\sigma_2}\setminus \left\{0\right\}\right)\times (0,1)$, 
 \[
  \begin{split}
  \gamma(\sigma_2,\eps)
  % \le F_{1}\left(U; C_{1/t_1}\right) &= \frac{1}{2}\int_{C_{1/t_2}} \abs{\nabla_\eps U_2}^2\d x + \frac{1}{4}\int_{B_{1/t_2}}\left(1- \abs{\Pi(U_2)}^2\right)\d x' \\
  % &\qquad  + \frac{1}{2}\int_{C_{1/t_1}\setminus C_{1/t_2}}\abs{\nabla_\eps H(\cdot,0)}^2 \d x\\
  \le F_{\eps}\left(U; C_{\sigma_2}\right) 
  &= F_{\eps}\left(U_1; C_{\sigma_1}\right) 
   + \frac{1}{2}\int_{C_{\sigma_2}\setminus C_{\sigma_1}}
   \abs{\nabla_\eps H}^2 \d x\\
  &= \gamma(\sigma_1, \eps) + \frac{1}{2}\int_{C_{\sigma_2}\setminus C_{\sigma_1}}\abs{\nabla_\eps H}^2 \d x.
  \end{split}
 \]
 The map $H$ does not depend on the $x_3$-variable and, as noticed, $H\equiv \Pi(H)$ on  $\left(B_{\sigma_2}\setminus \left\{0\right\}\right)\times (0,1)$. Thus
 \[
  \int_{C_{\sigma_2}\setminus C_{\sigma_1}}\abs{\nabla_\eps H}^2\d x = 
  \int_{B_{\sigma_2}\setminus B_{\sigma_1}}\abs{D\left(\dfrac{x'}{\abs{x'}}\right)}^2\d x' =2\pi \log\left(\dfrac{\sigma_2}{\sigma_1}\right) \! . 
 \]
 Therefore, summing up, we have that 
 \begin{equation}
  \label{eq:upperb1}
  \gamma(\sigma_2,\eps) \le \gamma(\sigma_1,\eps) 
   + \pi \log\left(\dfrac{\sigma_2}{\sigma_1}\right) \! ,
 \end{equation}
 which readily implies~$\tilde{\gamma}(\sigma_2,\eps) \leq \tilde{\gamma}(\sigma_1, \, \eps)$. 
\end{step}

\begin{step}[$\tilde{\gamma}$ is bounded]
 The limit in~\eqref{eq:conv_core} exists by monotonicity. It remains to prove that the limit is finite. As in the Ginzburg-Landau case, treated in \cite{BBH}, this follows from uniform lower and upper
 bounds for~$\tilde{\gamma}$.
 % The details are as follows. 
 Indeed, we fix~$\eps > 0$, $\sigma > 0$ %with~$0 < \eps < 1$
 and we let $U = U_\eps$ be the minimizer for $\gamma(\sigma,\eps)$. 
 Since (see e.g.~\cite[Theorem V.3]{BBH}, combined with a scaling argument)
 \[
  GL_\eps\left(\bar{u}_\eps;B_\sigma\right)\ge \pi \log\frac{\sigma}{\eps} - C,
 \]
 the inequality~\eqref{eq:GLbound2}, %and the assumption~\eqref{eq:ipeta}, 
 for~$\eps$ sufficiently small, implies
 \[
  \gamma(\sigma,\eps)
  = F_\eps\left(U_\eps;C_\sigma\right) \ge \pi \abs{\log\eps}-C,
 \]
 which implies that~$\tilde{\gamma}$ is bounded from below.
 As for the upper bound, let~$\zeta\in H^1(B_1;\S^2)$ be 
 any map such that~$\zeta(x^\prime) = (x^\prime, 0)$ 
 for~$x^\prime \in \partial B_1$.
 Given~$\eps>0$ and~$\sigma>0$, we define
 $V_\eps\in H^1(B_\sigma;\S^2)$ as
 \[
  V_\eps(x) :=
  \begin{cases}
   H(x) &\textrm{if } \eps \leq \abs{x^\prime} \leq \sigma \\
   \zeta\left(\dfrac{x^\prime}{\eps}\right) 
    &\textrm{if} \abs{x^\prime}\leq \eps
  \end{cases}
 \]
 Then, an explicit computation shows that
 \[
  \gamma(\sigma, \eps) \leq F_\eps(V_\eps;B_\sigma)
  = \pi\log\frac{\sigma}{\eps} 
   + \frac{1}{2}\int_{B^1} \abs{\nabla\zeta}^2 \d x^\prime
 \]
 which implies that~$\tilde{\gamma}$ is bounded from above.
\end{step}

\begin{step}[Conclusion]
 Given any sequence~$\eps_k\to 0$, by Helly's selection 
 theorem we can extract a (non-relabelled)
 subsequence in such a way that the sequence of
 functions~$\sigma\mapsto\tilde{\gamma}(\sigma,\eps_k)$ converge 
 pointwise to a limit~$\sigma\mapsto\tilde{\gamma}(\sigma,0)$.
 Moreover, the function~$\sigma\mapsto\tilde{\gamma}(\sigma,0)$ 
 is bounded and monotone, so the limit~$\gamma$ in~\eqref{eq:conv_core}
 exists in~$\R$. 
 \qedhere
\end{step}
\end{proof}

\begin{remark} \label{rk:linearcase}
 Suppose that~$\eta(\eps) = k\, \eps$,
 where~$k\in (0, \, 1/\sqrt{2}]$ is a constant 
 that does \emph{not} depend on~$\eps$.
 Then, a scaling argument shows that
 \[
  \gamma(\sigma,\eps) = \gamma\left(\frac{\sigma}{\eps},1\right) = \gamma\left(1,\frac{\eps}{\sigma}\right) =:
  \bar{\gamma}\left(\frac{\eps}{\sigma}\right)
 \]
 and, by Lemma~\ref{lemma:core}, the function
 $t\mapsto \bar{\gamma}(t) +\pi\log t$ is monotone nondecreasing
 and bounded. As a consequence, the limit
 \[
  \gamma = \lim_{\sigma\to 0}\lim_{\eps\to 0}
   \left(\gamma(\sigma,\eps) - \pi\log\frac{\sigma}{\eps}\right)
  = \lim_{t\to 0}
   \left(\bar{\gamma}\left(t\right) + \pi\log t\right)
 \]
 exists in~$\R$, with no need to pass to subsequences. 
\end{remark}

\subsection{An interpolation lemma}

This section is devoted to the following lemma,
which will be useful in our $\Gamma$-convergence analysis later on.
As before, we let~$C_\sigma := B_\sigma\times (0, 1)$ for~$\sigma > 0$.

\begin{lemma} \label{lemma:interp}
 Assume that~$\eta(\eps)\to 0$ as~$\eps\to 0$.
 For any~$M > 0$, $\sigma\in (0, 1)$ and~$\lambda\in (0, 1/2)$,
 there exists a number~$\eps_0 = \eps_0(M,\sigma,\lambda)> 0$
 such that the following statement holds.
 Let~$\eps\in (0, \eps_0)$, $h\in H^1(\partial B_\sigma; \S^1)$ 
 and~$U\in H^1(\partial B_\sigma\times (0, 1);\S^2)$
 be given. Suppose that the map
 $\bar{u} := \int_0^1 \Pi(U) \, \d x_3\colon\partial B_\sigma\to\R^2$
 is nowhere zero and has the same degree as~$h$ on~$\partial B_\sigma$, 
 and that
 \begin{equation} \label{hp:interp}
  \int_{\partial B_\sigma} \abs{\partial_\tau h}^2 \d s +
  \int_{\partial B_\sigma\times (0, 1)} \abs{\nabla_\eps U}^2 \d S
   + \frac{1}{\eps^2} \int_{\partial B_\sigma\times \{0,1\}}
   \abs{(U, \nu)}^2 \d s \leq M^2
 \end{equation}
 where~$\partial_\tau$ is the tangential derivative along~$\partial B_\sigma$. 
 Then, there exists~$W\in H^1(C_{\sigma(1 + \lambda)}\setminus C_{\sigma}; \S^2)$ such that
 \begin{align}
  &W(x) = U(x) \hspace{3cm} 
   \textrm{for } x\in \partial B_\sigma\times (0,1) \label{interp:bd1}\\
  &W(x) = \left(h\left(\frac{\sigma x'}{\abs{x'}}\right) \!, \, 0\right)\hspace{1.3cm}
   \textrm{for } x\in \partial B_{\sigma(1 + \lambda)}\times (0,1) 
   \label{interp:bd2} \\
  &F_\eps(W; C_{\sigma(1 + \lambda)}\setminus C_{\sigma}) 
  \leq \lambda\sigma M^2 + \frac{1}{\lambda\sigma} \int_{\partial B_\sigma} \abs{\bar{u} - h}^2 \d s \label{interp:energy}
 \end{align}
\end{lemma}
\begin{proof}
 Throughout the proof, we will use the notation~$A\lesssim B$
 as a short-hand for~$A \leq C B$, where~$C$ is a universal constant.
 Let
 \[
  a: = \sigma + \eta(\eps), \qquad b := a + \eps,
   \qquad c := \sigma(1 + \lambda).
 \]
 If~$\eta(\eps)\to 0$ as~$\eps\to 0$ and~$\eps$ is 
 small enough compared to~$\lambda\sigma$,
 then~$b < c$ and, in fact, $\lambda\sigma \lesssim c - b$.
 
 \setcounter{step}{0}
 \begin{step}[Construction of~$W$ on~$C_a\setminus C_\sigma$]
  First of all, by reasoning along the lines of Lemma~\ref{eq:lemub} 
  and applying~\eqref{hp:interp}, we have
  \begin{equation} \label{interp11}
   \begin{split}
    \frac{1}{2\eps^2} \int_{\partial B_\sigma} 
     \abs{(U(\cdot, x_3), \hat{e}_3)}^2 \d s 
    &\leq \frac{1}{\eps^2} \int_{\partial B_\sigma} 
     \abs{(U(\cdot,0), \hat{e}_3)}^2 \d s 
     + \frac{1}{\eps^2} \int_{\partial B_\sigma} 
     \abs{U(\cdot, x_3) - U(\cdot,0)}^2 \d s\\
    &\lesssim \frac{1}{\eps^2} \int_{\partial B_\sigma} 
     \abs{(U(\cdot,0), \hat{e}_3)}^2 \d s 
     + \int_{\partial B_\sigma\times (0, 1)} 
     \abs{\frac{\partial U}{\partial x_3}}^2 \d S \lesssim M^2
   \end{split}
  \end{equation}
  for a.e.~$x_3\in (0, \, 1)$. Next, we select a 
  value~$\bar{x}_3\in (0, \, 1)$ such 
  that~$U(\cdot,\bar{x}_3)\in H^1(\partial B_\sigma;\S^2)$ and
  \begin{equation} \label{interp12}
   \int_{\partial B_\sigma} \abs{\nabla_\eps U(\cdot,\bar{x}_3)}^2 \d s
   \leq \int_{\partial B_\sigma\times (0, \, 1)} 
    \abs{\nabla_\eps U(\cdot,\bar{x}_3)}^2 \d S
   \leq M^2
  \end{equation}
  Such a~$\bar{x}_3$ exists, by the
  assumption~\eqref{hp:interp} and Fubini's theorem.
  Now, let
  \[
   y(x) := x_3 - \sign(x_3 - \bar{x}_3)\frac{\abs{x'} - \sigma}{\eta(\eps)}
  \]
  For each given value of~$x_3$, as the value of~$\abs{x'}$
  increases in the range~$\sigma \leq \abs{x'}
  \leq \sigma + \eta(\eps)\abs{x_3 - \bar{x}_3}$,
  the function~$y(x)$ interpolates linearly between~$x_3$ and~$\bar{x}_3$. 
  We define~$W\colon C_a\setminus C_\sigma\to\S^2$ as follows:
  \[
   W(x) := 
   \begin{cases}
    U\left(\dfrac{\sigma x'}{|x'|}, y(x)\right) 
     &\textrm{if } \sigma \leq \abs{x'}
  \leq \sigma + \eta(\eps)\abs{x_3 - \bar{x}_3} \\[9pt]
    U\left(\dfrac{\sigma x'}{|x'|}, \bar{x}_3\right) 
     &\textrm{otherwise} 
   \end{cases}
  \]
  By construction, the function~$W$ belong to~$H^1(C_a\setminus C_\sigma;\S^2)$
  (in particular, its traces on either side of the 
  surface~$\abs{x'} = \sigma + \eta(\eps)\abs{x_3 - \bar{x}_3}$ 
  agree with one another).
  Moreover, $W$ satisfies~\eqref{interp:bd1} and
  \begin{equation} \label{interpWb}
   W(x) = W_a(x') := U\left(\frac{\sigma x'}{\abs{x'}}, \bar{x}_3\right)
   \qquad \textrm{for } \abs{x'} = a = \sigma + \eta(\eps).
  \end{equation}
  An explicit computation shows that
  \[
   \abs{\nabla_\eps W(x)} \lesssim 
   \abs{\nabla_\eps U \left(\frac{\sigma x'}{\abs{x'}}, x_3\right)}
    + \abs{\nabla_\eps U\left(\frac{\sigma x'}{\abs{x'}}, \bar{x}_3\right)}
  \]
  Therefore, keeping~\eqref{interp11} and~\eqref{interp12} in mind,
  we obtain the estimate
  \begin{equation} \label{interp1}
   F_\eps(W;C_a\setminus C_\sigma) \lesssim \eta(\eps) M^2.
  \end{equation}
 \end{step}
 
 \begin{step}[Construction of~$W$ on~$C_b\setminus C_c$]
  In view of~\eqref{interp11} and~\eqref{interp12},
  the function~$W_a$ defined in~\eqref{interpWb}
  is independent of the~$x_3$-variable and satisfies
  \begin{equation} \label{interp21}
   \int_{\partial B_a} \abs{\nabla_\eps W_a}^2 \d s
     + \frac{1}{\eps^2} \int_{\partial B_a} 
     \abs{(W_a, \hat{e}_3)}^2 \d s \lesssim M^2
  \end{equation}
  We recall the interpolation inequality
  \begin{equation} \label{interpGN}
   \norm{f}_{L^\infty(\partial B_a)}
   \lesssim \norm{f}_{L^2(\partial B_a)}^{1/2} 
    \norm{\partial_\tau f}_{L^2(\partial B_a)}^{1/2} 
    + a^{-1/2} \norm{f}_{L^2(\partial B_a)}
  \end{equation}
  which is valid\footnote{For the reader's convenience,
  we recall the proof of~\eqref{interpGN}. By a scaling argument,
  we can assume without loss of generality that~$a = 1$.
  By Sobolev embedding, we have
  \[
   \begin{split}
    \norm{f}_{L^\infty(\partial B_1)}^2
   = \norm{f^2}_{L^\infty(\partial B_1)}
   \lesssim \norm{\partial_\tau (f^2)}_{L^1(\partial B_1)}  
    + \norm{f^2}_{L^1(\partial B_1)}
   \lesssim \norm{f}_{L^2(\partial B_1)} 
    \norm{\partial_\tau f}_{L^2(\partial B_1)}  
    + \norm{f}_{L^2(\partial B_1)}^2
   \end{split}
  \]
  which implies~\eqref{interpGN}.}
  for any~$a > 0$ and~$f\in H^1(\partial B_b; \R)$.
%   here~$\partial_\tau f$ is the tangential derivative of~$f$
%   along~$\partial B_b$. 
  Due to~\eqref{interp21} and~\eqref{interpGN}, we have
  \begin{equation} \label{interp22}
   \norm{(W_a, \hat{e}_3)}_{L^\infty(\partial B_a)} 
    \lesssim \eps^{1/2} \left(1 + \left(\frac{\eps}{a}\right)^{1/2}\right) M 
  \end{equation}
  In particular, if~~$\eps$ is small enough 
  with respect to~$M$ and~$\sigma$, we have 
  ${(W_a, \hat{e}_3)} \leq 1/\sqrt{2}$
  everywhere in~$\partial B_a$, hence $\abs{\Pi(W_a)} \geq 1/\sqrt{2}$.
  We define
  \begin{equation} \label{interpWc}
   W_b(x') := \left(\frac{\Pi(W_a(x'))}{\abs{\Pi(W_a(x'))}}, \, 0\right)
   \qquad \textrm{for } x'\in \R^2\setminus\{0\},
  \end{equation}
  The map~$W_b$ is homogeneous
  (i.e., $W_b(s x') = W_b(x')$ for any~$s> 0$ 
  and~$x'\in\R^2\setminus\{0\}$), because so is~$W_a$.
  We also define $t(x') := \frac{b - \abs{x'}}{b - a}$ and 
  \[
   W(x) := \frac{t(x') \, W_a(x') + (1 - t(x')) \,  W_b(x')}
   {\abs{t(x') \, W_a(x') + (1 - t(x')) \,  W_b(x')}} 
   \qquad \textrm{for } x\in C_b\setminus C_a.
  \]
  This map is well-defined, thanks to~\eqref{interp22}. 
  It agrees with~$W_a$ on~$\partial B_a\times (0, 1)$,
  and with~$W_b$ on~$\partial B_b\times (0, 1)$.
  Moreover, recalling that~$b = a + \eps$, we have
  \[
   \abs{\nabla_\eps W} \lesssim  \abs{\nabla_\eps W_a}
   + \frac{1}{\eps} \abs{W_a - W_b}
   \lesssim \abs{\nabla_\eps W_a}
   + \frac{1}{\eps} \abs{(W_a, \hat{e}_3)}
  \]
  Then, with the help of~\eqref{interp21}, we conclude that
  \begin{equation} \label{interp2}
   F_\eps(W;B_b\setminus B_a) \lesssim \eps M^2.
  \end{equation}
 \end{step}
 
 \begin{step}[Construction of~$W$ on~$B_c\setminus B_b$]
  We observe that the maps $\bar{u} := \int_0^1 \Pi(U) \, \d x_3\colon\partial B_\sigma\to\R^2$, $h\colon \partial B_\sigma\to\S^1$ and~$\Pi(W_b)\colon\partial B_\sigma\to\S^1$ are all continuous, by Sobolev embeddings.
  Since we have assumed that~$\bar{u}$ has the same degree 
  as~$h$, and since~$\Pi(W_b)$ is homotopic to~$\bar{u}$
  by construction, lifting results imply that we can write
  \[
   \Pi(W_b(x')) = \left( e^{i\alpha(x')} h\left(\frac{x'}{|x'|}\right), \, 0\right)
   \qquad \textrm{for } x'\in\partial B_b,
  \]
  for some (real-valued) function $\alpha\in H^1(\partial B_b)$.
  Let~$\beta\in H^1(B_c\setminus B_b)$ be 
  the unique harmonic function that satisfies $\beta = 0$
  on~$\partial B_c$, $\beta = \alpha$ on~$\partial B_b$.
  We define
  \[
   W(x) := \left( e^{i\beta(x')} h\left(\frac{x'}{|x'|}\right), \, 0\right)
   \qquad \textrm{for } x\in C_c\setminus C_b.
  \]
  The map~$W$ agrees with~$W_b$ on~$\partial B_b\times (0, 1)$,
  satisfies~\eqref{interp:bd2} and
  \begin{equation*}
   \begin{split}
    F_\eps(W; C_c\setminus C_b)
    = \frac{1}{2} \int_{B_c\setminus B_b}
     \abs{D\Pi(W)}^2 \d x'
    &\lesssim \lambda\sigma\int_{\partial B_\sigma} \abs{\partial_\tau h}^2\d s
    + \int_{B_c\setminus B_b} \abs{D \beta}^2 \d x'
   \end{split}
  \end{equation*}
  By construction, $\beta$ is a minimizer of 
  the Dirichlet energy~$f\mapsto \int_{B_c\setminus B_b}\abs{D f}^2 \d x'$
  subject to its own boundary condition. Therefore,
  the integral of~$\abs{D\beta}^2$ can be estimated from above
  by considering a suitable competitor ---
  for instance, the function obtained 
  by interpolating linearly between~$\alpha$ and~$0$
  along the radial direction. Thus, we see that
  \begin{equation*}
   \begin{split}
    F_\eps(W; C_c\setminus C_b)
    &\lesssim \lambda\sigma \int_{\partial B_\sigma} 
     \left(\abs{\partial_\tau h}^2 + \abs{\partial_\tau \Pi(W_b)}^2\right)\d s
    + \frac{1}{\lambda\sigma} 
     \int_{\partial B_\sigma} \abs{h - \Pi(W_b)}^2 \d s 
   \end{split}
  \end{equation*}
  Then, keeping~\eqref{hp:interp} 
  and~\eqref{interp21} in mind, we can write
  \begin{equation*}
   \begin{split}
    F_\eps(W; C_c\setminus C_b)
    &\lesssim \lambda\sigma M^2
     + \frac{1}{\lambda\sigma} \int_{\partial B_c} 
     \left(\abs{h - \bar{u}}^2 + \abs{\bar{u} - \Pi(W_a)}^2
     + \abs{W_a - W_b}^2 \right) \d s  
   \end{split}
  \end{equation*}
  The integral of~$\abs{\bar{u} - \Pi(W_a)}^2$ at 
  the right-hand side can be estimated further by reasoning 
  as in Lemma~\ref{eq:lemub}, while the integral of~$\abs{W_a - W_b}^2$
  is bounded by~\eqref{interp21}. Thus, we obtain
  \begin{equation} \label{interp3}
   \begin{split}
    F_\eps(W; C_c\setminus C_b)
    &\lesssim \lambda\sigma M^2
     + \frac{1}{\lambda\sigma} \int_{\partial B_c} \abs{h - \bar{u}}^2\d s  
     + \frac{\eta(\eps)^2}{\lambda\sigma} M^2 
     + \frac{\eps^2}{\lambda\sigma} M^2
   \end{split}
  \end{equation}
  If we choose~$\eps$ small enough with respect to~$\lambda\sigma$, then~\eqref{interp:energy} follows from~\eqref{interp1}, 
  \eqref{interp2} and~\eqref{interp3}.
  \qedhere
 \end{step}
\end{proof}

\subsection{$\Gamma$-convergence analysis}
\label{ssec:Gamma}

The aim of this section is to prove a $\Gamma$-convergence
result for the functional~$F_\eps$, under the assumption~\eqref{eq:ipeta}.
Similar results are known in the context of the 
Ginzburg-Landau functional (see e.g.~\cite{AlicandroPonsiglione}).
First of all, we state and prove a compactness result.

\begin{theorem}\label{th:compactness}
Assume that the boundary datum~$G$ takes the form~\eqref{eq:bddatum},
with~$g\in H^{1/2}(\partial\Omega;\S^1)$, 
and that~\eqref{eq:ipeta} holds.
Let $M\in \mathbb{N}$ and $U_\eps\in\mathcal{A}_G$
be a sequence that satisfies 
\begin{equation}
\label{eq:energy_bound}
F_\eps(U_\eps) \le M\pi \abs{\log\eps} + C
\end{equation}
for some~$\eps$-independent constant~$C$.
Then, the following properties hold.
\begin{enumerate}[label=(\roman*)]
 \item We have
  \begin{equation*}
   \begin{split}
    &\frac{\partial U_\eps}{\partial x_3}\xrightarrow{\eps\to 0}0 \qquad \textrm{strongly in } L^2(Q;\mathbb{R}^3)\label{eq:not_tangent}\\
    &U_{\perp,\eps}\xrightarrow{\eps\to 0} 0 \qquad \textrm{strongly in } L^2(Q;\R).
  \end{split}
\end{equation*}

 \item There exist a subsequence of $\eps$,
 an $N$-uple of distinct points $a_1,\ldots,a_N$ in $\bar{\Omega}$ 
 and nonzero integers $d_1,\ldots, d_N$ 
 such that $T:= \sum_{n=1}^N \abs{d_n} \le M$, 
 \begin{equation}
 \label{eq:poin-hopf}
  \sum_{n=1}^{N} d_n = \deg(g,\partial\Omega),
 \end{equation}
 and
 \begin{equation}
 \label{eq:vorticity_convergence}
 J \bar{u}_\eps \rightarrow 2\pi\sum_{n=1}^N d_n \delta_{a_n} 
 \qquad \textrm{strongly in } \left(W^{1,\infty}(\Omega)\right)'
 \textrm{ as } \eps\to 0,
 \end{equation}
 where $\bar{u}_\eps = \int_{0}^{1}\Pi(U_\eps)\d t$.

 \item Assume in addition that~$T = M$.
 Then, $N = T =M$ and, for any index~$n= 1, \ldots, N$,
 we have $\abs{d_n} = 1$ and $a_j\in \Omega$.
 Moreover, there exists a further (non-relabelled) subsequence
 and a map~$u$ that satisfies
 \begin{gather}
  u\in W^{1,p}(\Omega;\mathbb{S}^1) \quad \textrm{for any } p\in (1, \, 2) \textrm{ and } u\in H^1_{\mathrm{loc}}\left(\Omega\setminus \bigcup_{n=1}^N \left\{a_n\right\};\mathbb{S}^1\right) \label{eq:limitmap0} \\
  \bar{u}_\eps \rightharpoonup u \ \textrm{ weakly in } W^{1,p}(\Omega, \, \R^2) \textrm{ for any } p\in (1, \, 2) \textrm{ and in }
  H^1_{\mathrm{loc}}\left(\Omega\setminus \bigcup_{n=1}^N \left\{a_n\right\};\R^2\right) \label{eq:limitmap}
 \end{gather}
 and
 \begin{equation}
 \label{eq:limitdegree}
 \deg(u,\partial B_\sigma(a_n)) = d_n
 \end{equation}
 for any~$n=1, \ldots, N$ and any~$\sigma$ small enough.
 \end{enumerate}
\end{theorem}

\begin{proof}%[Proof of Theorem~\ref{th:compactness}]
 Let $U_\eps\in \mathcal{A}_{G}$ satisfy~\eqref{eq:energy_bound}.
% Since~$\int_{Q}\abs{U_\eps}^2 \d x = \abs{\Omega}$ 
% % (recall that $\abs{\Omega}=1$),
% there exists a (non-relabelled) subsequence of~$\eps\to 0$ 
% and a map $U^*\in L^2(Q)$ such that 
% \begin{equation}
% \label{eq:conv_1}
% U_\eps \weak U \quad \textrm{ weakly in } L^2(Q).
% \end{equation} 
 Then, keeping~\eqref{eq:ipeta} into account, we obtain
 \[
  \begin{split}
   &\frac{1}{2}\int_{Q}\abs{\frac{\partial U_\eps}{\partial x_3}}^2\d x  \le \eta^2(\eps)F_\eps(U_\eps)\le \eta^2(\eps)(M\pi \abs{\log\eps} + C) \to 0\\
   &\int_{\Omega\times \left\{0,1\right\}}\abs{(U_\eps,\nu)}^2 \d x'\le \eps^2 F_\eps(U_\eps)\le \eps^2(M\pi \abs{\log\eps} + C) \to 0
  \end{split}
 \]
 as~$\eps\to 0$.
 Therefore, thanks to~\eqref{eq:ipeta}, we immediately have
 \[
  \begin{split}
   &\frac{\partial U_\eps}{\partial x_3} \xrightarrow{\eps \to 0} 0 \quad\textrm{ strongly in } L^2(Q)\\
   &U_\eps - \Pi(U_\eps)\xrightarrow{\eps \to 0} 0 \quad \textrm{ strongly in } L^2(\Omega\times\left\{0,1\right\})
  \end{split}
 \]
 and Statement~(i) follows, because of the embeddings
 $H^1(0, 1; L^2(\Omega))\hookrightarrow C([0, \, 1]; \, L^2(\Omega))
 \hookrightarrow L^2(Q)$.
 Now, let  $\bar{u}_\eps:= \int_{0}^1\Pi(U_\eps)\,\d x_3$. 
 Due to~\eqref{eq:GLbound2} and~\eqref{eq:energy_bound}, we have
 \begin{equation} \label{compactness1}
  GL_\eps(\bar{u}_\eps) 
  % = \frac{1}{2}\int_{\Omega}\abs{D \bar{u}_\eps}^2 \d x' + \frac{1}{4\eps^2}\int_{\Omega}\left(1-\abs{\bar{u}_\eps}^2\right)^{\color{blue}2} \d x'
  \le M\pi\abs{\log\eps} + C.
 \end{equation}
 Thus, thanks to classical compactness results for the Ginzburg-Landau functional (see e.g. \cite[Theorem 6.1]{AlicandroPonsiglione}), 
 we have the desired conclusion, i.e. Statement~(ii).
 %\eqref{eq:poin-hopf} and \eqref{eq:vorticity_convergence}
 
 Similarly, the proof of Statement~(iii) follows by pre-existing 
 results on the Ginzburg-Landau functional, which apply
 to our setting because of~\eqref{compactness1}.
 In particular, assuming~$T = M$, a proof of
 the fact that~$T = N$ and~$\abs{d_n} = 1$,
 $x_n\in\Omega$ for any~$n$ can be found
 in~\cite[Theorem~5.3]{AlicandroPonsiglione}.
 The existence of~$u$ satisfying~\eqref{eq:limitmap0},
 \eqref{eq:limitmap} and~\eqref{eq:limitdegree}
 follows by energy estimates~\cite{Jerrard, Sandier};
 the details of the argument can be found, e.g.,
 in~\cite[Lemma~4.3 and~4.4]{AG-Reno}
%  (for a proof of the weak
%  convergence in~$H^1_{\mathrm{loc}}$) 
 and in~\cite[proof of Theorem~1.1, p.~1620]{struwe94}
 (for a proof of the weak convergence in~$W^{1,p}$).
\end{proof}

\begin{remark} \label{rk:compactness}
 Under the assumptions of Theorem~\ref{th:compactness},
 when~$T=M$ there exists a (nonrelabelled) subsequence such that
 \[
  U_\eps\to U \qquad \textrm{strongly in } L^q(Q) 
  \ \textrm{ for any } q < +\infty,
 \]
 where~$U(x) := (u(x^\prime), 0)$ for~$x=(x^\prime, x_3)\in Q$
 and~$u$ is given by Theorem~\ref{th:compactness}.
 Indeed, Theorem~\ref{th:compactness} and Lemma~\ref{eq:lemub},
 combined, imply that~$U_\eps\to U$ strongly in~$L^2(Q)$.
 The convergence in~$L^q(Q)$ for any finite~$q$ follows 
 by interpolation, as~$U_\eps$ takes values in~$\S^2$ and, hence,
 is uniformly bounded in~$L^\infty(Q)$.
%  {\BBB Non so se si possa arrivare facilmente ad avere convergenza
%  in $W^{1,p}$ o in~$H^1_{\mathrm{loc}}$, perch\'e il lemma~\ref{eq:lemub}
%  dà solo una stima in~$L^2$\ldots}
% %  {\color{red} Mettere che $U_\eps$ converge debole in $H^1_{loc}(Q\setminus \bigcup_{n=1}^N (B_\sigma(a_n)\times(0,1))$ e che $U= \bar{u}$ (vedi \eqref{eq:limitmap}), grazie a Lemma \ref{eq:lemub}}
\end{remark}

%\begin{step}{Preliminary computation.}

%For any $x_3\in [0,1]$
%\[
%\begin{split}
%\abs{\abs{\bar{u}(x')}^2-\abs{u(x',x_3)}^2} &\le 2\abs{\int_{0}^1\left(u(x',x_3)-u(x', y)\right)\d y}
%=2\abs{\int_{0}^1\left(\int_{y}^{x_3}\frac{\partial u}{\partial s}(x',s)\d s\right)\d y}\\
%&\le 2\int_{0}^1\abs{\frac{\partial u}{\partial s}(x',s)}\d s.
%\end{split}
%\]
%Therefore (see \eqref{eq:bound1})
%\[
%\int_{\Omega}\abs{\abs{\bar{u}(x')}^2-\abs{u(x',0)}^2} \d x'\le 2\left(\int_{Q}\abs{\frac{\partial u}{\partial x_3}}^2\d x\right)^{1/2}\le 2\eps F_\eps(U) 
%\]
%and thus
%\[
%\int_{\Omega}\left(1-\abs{\bar{u}}^2\right) \d x'  \le
%\int_{\Omega}\abs{\abs{\bar{u}(x')}^2-\abs{u(x',0)}^2} \d x' + \int_{\Omega}\left(1-\abs{u(x',0)}^2\right)\d x' \le 
%3\eps F_\eps(U)
%\]
%
%The Jensen's Inequality gives that 
%\begin{equation}
%\label{eq:Jenses}
%F_\eps(U)\ge \int_{Q}\abs{D u}^2\d x = \int_{\Omega}\left(\int_{0}^1 \abs{D u}^2\d x_3\right)\d x'\ge \int_{\Omega}\abs{D \bar u}^2 \d x'.
%\end{equation}
%Therefore from $F_\eps(U)\le C\abs{\log\eps}$  we deduce that 
%\begin{equation}
%\label{eq:GLbound}
%\frac{1}{2}\int_{\Omega}\abs{D \bar u}^2 \d x' + \frac{1}{\eps}\int_{\Omega}\left(1-\abs{\bar{u}}^2\right) \d x'\le C\abs{\log\eps}
%\end{equation}
%\end{step}

Finally, we state a compactness result,
in the sense of $\Gamma$-convergence, for our sequence
of functionals~$(F_\eps)_{\eps > 0}$.

\begin{theorem} \label{th:Gamma}
Assume that the boundary datum~$G$ takes the form~\eqref{eq:bddatum},
with~$g\in H^{1/2}(\partial\Omega;\S^1)$, 
and that~\eqref{eq:ipeta} holds. Then, 
for any sequence~$\eps_k\to 0$, there exists a 
(non-relabelled) subsequence and a number~$\gamma\in\R$
such that the following statements hold.

\begin{enumerate}[label=(\roman*), leftmargin=*]
\item \emph{$\Gamma$-liminf Inequality.}
Let~$\a = (a_1, \ldots, a_N)$ be an $N$-uple of distinct 
points in~$\Omega$, $\db = (d_1, \ldots, d_N)$ an~$N$-uple
of integers with~$\abs{d_n}=1$ for any~$n$, $\sum_{j=1}^n d_n = \deg(g,\partial\Omega)$, and let $U_{\eps_k}\in \mathcal{A}_{G}$ 
be a sequence that %for which $\bar{u}_{\eps_k} = \int_{0}^1\Pi(U_{\eps_k})$
satisfies~\eqref{eq:vorticity_convergence}. Then, 
\begin{equation}
\begin{split}
\label{eq:liminf}
\liminf_{k\to+\infty} \left(F_{\eps_k}(U_{\eps_k})- N\pi\abs{\log\eps_k}\right)
\ge W_g(\a, \db) + N\gamma. 
\end{split}
\end{equation}

\smallskip
\item \emph{$\Gamma$-limsup Inequality.}
Conversely, given an $N$-uple~$\a = (a_1, \ldots, a_N)$ of distinct 
points in~$\Omega$, an~$N$-uple~$\db = (d_1, \ldots, d_N)$
of integers with~$\abs{d_n}=1$ for any~$n$
and~$\sum_{j=1}^n d_n = \deg(g,\partial\Omega)$,
there exists a sequence~$U_{\eps_k}\in\mathcal{A}_G$ that satisfies \eqref{eq:vorticity_convergence} and
\begin{equation}
\label{eq:limsup}
\limsup_{k\to+\infty} \left(F_{\eps_k}(U_{\eps_k})- N\pi\abs{\log\eps_k}\right)
\le  W_g(\a, \db) + N\gamma.
\end{equation} 
\end{enumerate}
\end{theorem}

The number~$\gamma$ is the core energy,
as defined by Lemma~\ref{lemma:core}. As we have not investigated
the uniqueness of the limit in~\eqref{eq:core}, we cannot
exclude the possibility that~$\gamma$ --- and hence,
the~$\Gamma$-limit of the functionals~$F_{\eps_k}$
--- depends on the subsequence~$\eps_k$.
In order to identify the~$\Gamma$-limit, a 
more thorough analysis of the minimization 
problem~\eqref{eq:core} is in order. Nevertheless, 
Theorem~\ref{th:Gamma} is enough to complete the proof 
of Theorem~\ref{th:minimizers}, because the $\Gamma$-limits 
along two different subsequences will differ only for an 
additive constant, which does not affect minimizers of
the Renormalized Energy.

As already observed, if we tune $\eta(\eps)$ in such a way 
that $\eta(\eps) = k\,\eps$ for some constant~$k\in (0, \frac{1}{\sqrt{2}}]$,
then we have a $\Gamma$-convergence result in the usual sense --- namely, the lower and the upper bound hold for the entire sequence~$\eps\to 0$ and we uniquely identify the core energy~$\gamma$. 
See Remark~\ref{rk:linearcase} for the details.

%{\BBB \textbf{Si pu\`o osservare che, nel caso~$\eta(\eps) = k\eps$
%con~$k\in (0, \, \frac{1}{\sqrt{2}}]$ costante, 
%il Remark~\ref{rk:linearcase} consente di avre un risultato 
%di~$\gamma$-convergenza vero e proprio.}}

The proof of Theorem~\ref{th:Gamma} follows along the 
lines of~\cite[Theorem~5.3]{AlicandroPonsiglione}.
We recall a few facts, which will be useful in the proof.
Given~$0 < r < R$, consider the minimization problem
\[
 m := \min\left\{\frac{1}{2} \int_{B_R\setminus B_r} \abs{D v}^2 \d x'\colon
 v\in H^1(B_R\setminus B_r; \S^1), \ \deg(v, \partial B_R) = 1\right\} 
\]
An argument based on Jensen's inequality
(see e.g.~\cite[Lemma~1.1]{Sandier}) shows that 
$m = \pi\log(R/r)$, and the minimzers are 
all and only the maps $h_\alpha\colon B_R\setminus B_r\to\S^1\subset\C$
given by
\begin{equation} \label{h_alpha}
 h_\alpha(\rho e^{i\theta}) := \alpha e^{i\theta},
 \qquad \textrm{for } \rho > 0, \ 0\leq \theta \leq 2\pi,
\end{equation}
where the parameter~$\alpha\in\C$ is such that~$\abs{\alpha}=1$.

\begin{lemma} \label{lemma:delta}
 For any~$\delta > 0$, there exists a number~$c(\delta) > 0$
 such that the following property holds: for any~$R > 0$
 and any map~$v\in H^1(B_R\setminus B_{R/2};\S^1)$
 such that 
 \[
  \inf_{\alpha\in\S^1} \left(
   \norm{D v - D h_\alpha}_{H^1(B_{R}\setminus B_{R/2})}
   + R^{-1}\norm{v - h_\alpha}_{L^2(B_{R}\setminus B_{R/2})}\right) 
   \geq \delta,
 \]
 there holds
 \[
  \frac{1}{2} \int_{B_{R}\setminus B_{R/2}} \abs{D v}^2 \d x^\prime
  \geq \pi\log 2 + c(\delta)
 \]
\end{lemma}
\begin{proof}
 By a scaling argument, we can assume without loss of generality
 that~$R=1$. Then, the lemma follows by arguing by contradiction
 (see e.g.~\cite[proof of Theorem~5.3, (ii)]{AlicandroPonsiglione}).
\end{proof}

\begin{proof}[Proof of Theorem~\ref{th:Gamma}]
Let~$\eps_k\to 0$ be a given sequence. By Lemma~\ref{lemma:core},
we can extract a (non-relabelled) subsequence in such 
a way that the limit~$\gamma$, defined by~\eqref{eq:conv_core},
exists in~$\R$. We restrict our attention
to this particular subsequence~$(\eps_k)_{k\in\N}$
but we write~$\eps$ instead of~$\eps_k$, to simplify the notation.

\medskip
\setcounter{step}{0}
\begin{step}
{{\itshape Lower bound.}}
Let~$\a=(a_1, \ldots, a_N)$, $\db=(d_1, \ldots, d_N)$
and~$U_\eps\in\mathcal{A}_G$ satisfy the hypothesis.
We can certainly assume that $\Lambda := \liminf_{\eps\to 0} (F_\eps(U_\eps)
 - \pi N\abs{\log\eps})$ is finite,
otherwise there is nothing to prove.
Moreover, up to extraction of a subsequence, we can assume
\begin{equation} \label{eq:liminf0}
 \Lambda = \lim_{\eps\to 0}
 \left( F_\eps(U_\eps) - \pi N\abs{\log\eps} \right) < +\infty.
\end{equation}
This allows us to extract further subsequences, without changing 
the limit of the energies.
Let~$\bar{u}_\eps := \int_0^1\Pi(U_\eps) \, \d x_3$.
By Theorem~\ref{th:compactness}, we can extract a subsequence
in such a way that $\bar{u}_\eps$ converge to a limit 
$u\colon\Omega\setminus\{a_1, \ldots, a_N\}\to\S^1$,
which satisfies~\eqref{eq:limitmap0}, \eqref{eq:limitmap}
and~\eqref{eq:limitdegree}.
Take a number~$\sigma > 0$ be small enough, so that the balls~$B_\sigma(a_n)$
are pairwise disjoint and contained in~$\Omega$,
and let~$\Omega_\sigma := \Omega\setminus\cup_{n=1}^N B_\sigma(a_n)$.
Then, \eqref{eq:GLbound2}, \eqref{eq:limitmap} and~\cite[Theorem~I.7]{BBH} 
imply
\begin{equation}
 \label{eq:liminf1}
 \begin{split}
  F_\eps\left(U_\eps; \Omega_\sigma\times (0, \, 1)\right)
  \geq GL_\eps(\bar{u}_\eps; \Omega_\sigma)
  &\geq \frac{1}{2} \int_{\Omega_\sigma} \abs{D u} ^2 \d x^\prime
   + \mathrm{o}_{\eps\to 0}(1) \\
  &\geq \pi N \abs{\log\sigma} + W_g(\a, \, \db) + \mathrm{o}_{\sigma\to 0}(1)
   + \mathrm{o}_{\eps\to 0}(1) 
 \end{split}
\end{equation}
For each index~$n$, let~$C_\sigma(a_n) := B_\sigma(a_n)\times (0, \, 1)$.
To conclude the proof of~\eqref{eq:liminf}, it suffices to show that
\begin{equation}
 \label{eq:liminf2}
 \begin{split}
  \liminf_{\eps\to 0}\big(F_\eps\left(U_\eps; C_\sigma(a_n)\right)
  - \pi\abs{\log\eps}\big) \geq \pi \log\sigma
  + \gamma + \mathrm{o}_{\sigma\to 0}(1)
 \end{split}
\end{equation}
for each~$n\in\{1, \ldots, N\}$.
For then, \eqref{eq:liminf}
would follow by combining~\eqref{eq:liminf1}
with~\eqref{eq:liminf2}, and taking the limit as~$\sigma\to 0$.

Let us fix and index~$n$ and assume, for simplicity, that~$d_n=1$.
(If~$d_n=-1$, a similar argument applies.)
%with~$\bar{U}_\eps := (U_{\eps,1}, -U_{\eps, 2}, U_{\eps, 3})$ in place of~$U_\eps$.)
Let~$\delta > 0$ be a small parameter, and let~$L\in\N$
be an integer. $\Gamma$-convergence results
for the Ginzburg-Landau functional 
(see e.g.~\cite[Theorem~5.3 and Remark~5.2]{AlicandroPonsiglione}) imply that 
\begin{equation} \label{eq:liminf3}
 GL_\eps(\bar{u}_\eps; B_{2^{-L}\sigma}(a_n))
 \geq \pi \log\frac{\sigma}{2^L \, \eps} - C_*
\end{equation}
where~$C_*$ is a universal constant
(in particular, independent of~$\eps$, $\sigma$, $L$).
Let us choose~$L = L(\delta)$ in such a way
that~$L c(\delta) \geq \gamma + C_*$,
where~$c(\delta)$ is given by Lemma~\ref{lemma:delta}.
For any~$\ell\in\{1,\ldots,L\}$, let~$\Gamma_\ell(a_n)
:= B_{2^{1-\ell}\sigma}(a_n)\setminus B_{2^{-\ell}\sigma}(a_n)\subset\R^2$.
Suppose first that
\begin{equation} \label{eq:liminf-cases}
 \inf_{\alpha\in\S^1} \left( 
  \norm{u - h_\alpha(\cdot - a_n)}_{L^2(\Gamma_\ell(a_n);\R^2)}
  + 2^{\ell - 1}\sigma^{-1} \norm{u
   - h_\alpha(\cdot - a_n)}_{L^2(\Gamma_\ell(a_n);\R^2)}\right)
 \geq \delta 
\end{equation}
for any~$\ell\in \{1,\ldots,L\}$,
where~$h_\alpha$ is defined in~\eqref{h_alpha}.
By applying~\eqref{eq:GLbound2} and~\eqref{eq:limitmap},
we deduce
\[
 \begin{split}
  F_\eps(U_\eps; C_\sigma(a_n))
  \geq GL_\eps(\bar{u}_\eps; B_\sigma(a_n))
  \geq GL_\eps(\bar{u}_\eps; B_{2^{-L}\sigma}(a_n))
  + \frac{1}{2} \sum_{\ell=1}^L \int_{\Gamma_\ell(a_n)} \abs{D u}^2 \d x^\prime
  + \mathrm{o}_{\eps\to 0}(1)
 \end{split}
\]
Then, \eqref{eq:liminf3} and Lemma~\ref{lemma:delta} imply
\[
 \begin{split}
  F_\eps(U_\eps; C_\sigma(a_n))
  &\geq \pi \log\frac{\sigma}{2^L \, \eps} 
   - C_* + L \left(\pi\log2 + c(\delta)\right)
   + \mathrm{o}_{\eps\to 0}(1) 
  \geq \pi \log\frac{\sigma}{\eps} +\gamma
   + \mathrm{o}_{\eps\to 0}(1),
 \end{split}
\]
because of our choice of~$L$. This proves~\eqref{eq:liminf2},
in case~\eqref{eq:liminf-cases} holds.

Next, we consider the case~\eqref{eq:liminf-cases} does not hold
--- which means, there exist an index~$\ell$ and a value~$\alpha\in\S^1$
such that 
\begin{equation} \label{eq:liminf3.5}
 \norm{u - h_\alpha(\cdot - a_n)}_{L^2(\Gamma_\ell(a_n);\R^2)}
  + 2^{\ell - 1}\sigma^{-1} \norm{u
   - h_\alpha(\cdot - a_n)}_{L^2(\Gamma_\ell(a_n);\R^2)}
 \leq \delta
\end{equation}
By applying~\eqref{eq:GLbound2} and~\eqref{eq:liminf0}, and
reasoning as in~\eqref{eq:liminf1}, \eqref{eq:liminf3}, we deduce
(here $\Omega_{2^{1 - \ell}\sigma} := \Omega\setminus\cup_{n=1}^N B_{2^{1-\ell}\sigma}(a_n)$):
\begin{equation} \label{eq:liminf4}
 \begin{split}
 F_\eps(U_\eps; \Gamma_\ell(a_n))
 &\leq \pi N \abs{\log\eps} + \Lambda + 1
 - GL_\eps(\bar{u}_\eps; \Omega_{2^{1 - \ell}\sigma})
 - \sum_{n=1}^N GL_\eps(\bar{u}_\eps; B_{2^{-\ell}\sigma}(a_n)) \\
 &\leq \pi N \abs{\log\eps} + \Lambda + 2
 - \pi N \abs{\log(2^{1 - \ell}\sigma)} 
 - \pi N \log\frac{\sigma}{2^{\ell}\eps} + N C_* \leq C.
 \end{split}
\end{equation}
Here and in what follows, $C$ denotes a generic constant
which may depend on~$\Lambda$, $N$ and~$(\a, \db)$,
but not on~$\eps$, $\sigma$, $\delta$.
Thanks to~\eqref{eq:liminf3.5}, \eqref{eq:liminf4}
and Fubini's theorem, we find~$\bar{\sigma} \in (2^{1-\ell}\sigma, 
\, 2^{-\ell}\sigma)$ such that
\[
 F_\eps(U_\eps; \partial B_{\bar{\sigma}}(a_n)) \leq \frac{C}{\bar{\sigma}},
 \qquad \int_{\partial B_{\bar{\sigma}}(a_n)} \abs{\bar{u}_\eps - h}^2 \d s \leq C\delta,
\]
where~$h:= h_\alpha(\cdot - a_n)$ and~$\alpha\in\S^1$
is as in~\eqref{eq:liminf-cases}. 
By applying Lemma~\ref{lemma:interp},
we find a map $W_\eps\in H^1(C_{\bar{\sigma}(1 + \delta^{1/2})}(a_n) \setminus C_{\bar{\sigma}}(a_n); \S^2)$ such that $W_\eps(x) = (h(x'), 0)$
if~$x \in \partial B_{\bar{\sigma}(1 + \delta^{1/2})}(a_n)\times (0,1)$,
$W_\eps = U_\eps$ on~$\partial B_{\bar{\sigma}}(a_n)\times 
(0,1)$ and
\begin{equation}
 \label{eq:liminf5}
 F_\eps(W_\eps; C_{ \bar{\sigma}(1 + \delta^{1/2})}(a_n) 
  \setminus C_{\bar{\sigma}}(a_n))
  \leq C \delta^{1/2} \left(1 + \frac{1}{\bar{\sigma}}\right)
\end{equation}
Let $R\colon\R^3\to\R^3$ be the rotation given
by~$R(x) := (\alpha^{-1} x', x_3)$ and let
$V_\eps\colon B_{\bar{\sigma}(1 + \delta^{1/2})}\to\S^2$
be given by
\[
 V_{\eps}(x) :=
 \begin{cases}
  R W_\eps(a_n + x) &\textrm{if } \bar{\sigma}
   < \abs{x'} \leq \bar{\sigma}(1 + \delta^{1/2}) \\
  R U_\eps(a_n + x) 
   &\textrm{if } \abs{x'} \bar{\sigma}.
 \end{cases}
\]
By construction, the map~$V_\eps$
satisfies~$V_\eps(x) = (\frac{x'}{\abs{x'}}, 0)$
for~$x\in\partial B_{\bar{\sigma}(1 + \delta^{1/2})}$. Therefore,
observing that the functional~$F_\eps$ is invariant
under isometries and recalling~\eqref{eq:liminf5}, we obtain
\begin{equation*}
 \begin{split}
  \gamma(\bar{\sigma}(1 + \delta^{1/2}), \eps)
  \leq F_\eps(V_\eps; \, C_{\bar{\sigma}(1 + \delta^{1/2})}(a_n))
  &\stackrel{\eqref{eq:liminf5}}{\leq} 
   F_\eps(U_\eps; \, C_{\bar{\sigma}}(a_n))
   + C \delta^{1/2} \left(1 + \frac{1}{\bar{\sigma}}\right) 
 \end{split}
\end{equation*}
where~$\gamma(\bar{\sigma}, \, (1 - \delta^{1/2})\eps)$
is defined by~\eqref{eq:core}.
Therefore, Lemma~\ref{lemma:core} gives
\begin{equation} 
 \label{eq:liminf-quasi}
 \begin{split}
  \liminf_{\eps\to 0}\big(F_\eps(U_\eps; \, C_{\bar{\sigma}}(a_n))
   - \pi\abs{\log\eps}\big)
  &\geq \pi\log(\bar{\sigma}(1 + \delta^{1/2})) + \gamma + \mathrm{o}_{\sigma\to 0}(1)
   - C \delta^{1/2} \left(1 + \frac{1}{\bar{\sigma}}\right)
 \end{split}
\end{equation}
Finally, the estimate~\eqref{eq:GLbound2}, combined with~\eqref{eq:limitmap}
and Jensen's inequality (see e.g.~\cite[Lemma~1.1]{Sandier}) imply
\begin{equation} 
 \label{eq:liminf-ancorapoco}
 \begin{split}
  \liminf_{\eps\to 0} F_\eps(U_\eps; \, C_\sigma(a_n) \setminus C_{\bar{\sigma}}(a_n))
  &\geq \liminf_{\eps\to 0} GL_\eps(\bar{u}_\eps; \, B_\sigma(a_n) \setminus B_{\bar{\sigma}}(a_n)) \\
  &\geq \frac{1}{2} \int_{B_\sigma(a_n) \setminus B_{\bar{\sigma}}(a_n)}
  \abs{D u}^2 \d x'
  \geq \pi\log\frac{\sigma}{\bar{\sigma}}
 \end{split}
\end{equation}
Combining~\eqref{eq:liminf-quasi} with~\eqref{eq:liminf-ancorapoco}, 
and taking the limit as~$\delta\to 0$ (but keeping~$\sigma$ fixed!),
\eqref{eq:liminf2} follows.

% Let $U_\eps\in \mathcal{A}_{G}$ be a sequence verifying the hypothesis. 
% We have that 
% \begin{equation}
% \label{eq:liminf1}
% F_\eps(U_\eps) -N\pi\abs{\log\eps} 
% \ge GL_\eps(\bar{u}_\eps) -N\pi\abs{\log\eps} + o_{\eps\to 0}(1)\ge W_g(a,d) + N\gamma +
%  o_{\eps\to 0}(1),
% \end{equation}
% thanks to \cite[Theorem 6.1]{AlicandroPonsiglione}.
\end{step}

\medskip
\begin{step}
{{\itshape Upper bound.}}
Suppose we are given $N$ distinct points $\a =(a_1,\ldots,a_N)$
in~$\Omega$ and integers $\db=(d_1,\ldots,d_N)$ 
with $\abs{d_n}=1$ for any $n$ and such 
that~$\sum_{n=1}^N d_n = \deg(g,\partial \Omega)$. 
To simplify the notation, we assume first that~$d_n = 1$ for any~$n$. 
% Without loss of generality we can assume that $d_n = 1$ for all~$n$. 
First of all we construct a sequence~$U_{\eps,\sigma}$ depending on~$\eps$ and on a parameter~$\sigma$ that, as in Lemma~\ref{lemma:core},
measures the the radii of the small balls around the vortices where the core energy concentrates. 
Then, a diagonal argument (see Lemma~\ref{lemma:recovery} below) will produce our recovery sequence. 

We construct the sequence $U_{\eps,\sigma}$ as follows. First of all we note that there is a unique canonical harmonic map~$u^*$
associated with~$(\a, \db, g)$ \cite[Section~I.3]{BBH}.
We let 
\begin{equation}
\label{eq:cyl-harmo-map}
U^*(x', x_3):= (u^*(x'), 0)\qquad \textrm{ for any } x=(x',x_3) \in \left(\Omega\setminus \bigcup_{n=1}^N\left\{a_n\right\}\right)\times (0,1).
\end{equation}
We take~$\sigma>0$ so small that the closed balls $\bar{B}_{\sigma}(a_n)$ are  mutually disjoint and contained in~$\Omega$ and we consider the perforated domain~$Q_\sigma$ defined by 
\[
Q_\sigma:= Q\setminus \bigcup_{n=1}^N C_{2\sigma}(a_n),
\qquad \textrm{where } 
C_{2\sigma}(a_n) := B_{2\sigma}(a_n)\times (0,1).
\]
Then, we let 
\begin{equation}
\label{eq:recovery_fuori}
U_{\eps,\sigma} = U^* \qquad \textrm{ on } Q_\sigma.
\end{equation}
Since $U^*$ does not depend on~$x_3$, 
$(U^*,\hat{e}_3)=0$ everywhere in~$x\in Q_\sigma$ and
$d_n = 1$ for any~$n$, we have that (see~\eqref{eq:RE_def})
\begin{equation}
\label{eq:energyU}
F_\eps(U^*;Q_\sigma) = \frac{1}{2}\int_{\Omega\setminus \bigcup_{n=1}^N B_{2\sigma}(a_n)}\abs{D u^*}^2\d x'
= \pi N\abs{\log(2\sigma)} + W_g(\a,\db) 
+ \mathrm{o}_{\sigma\to 0}(1)
\end{equation}

\medskip
Now we work in the sets $C_{2\sigma}(a_n)\setminus C_{\sigma}(a_n)$ for $n=1,\ldots,N$. 
For each fixed~$n$, %such that~$d_n = 1$, 
we let~$u_0$ be the hedgehog centered in $a_n$, namely 
\[
u_0(x') = \frac{x'-a_n}{\abs{x'-a_n}},\qquad 
\textrm{for } x'\in B_{2\sigma}(a_n)\setminus B_{\sigma}(a_n).
\]
% Instead, if~$d_n=-1$, we define~$u_0$ as
% \[
% u_0(x') = \frac{\overline{x'}-\overline{a_n}}{\abs{x'-a_n}},\qquad 
% \textrm{for } x'\in B_{\sqrt{\eps}}(a_n)\setminus B_{\sqrt{\eps}/2}(a_n).
% \]
% where the bars denote conjugation, in the sense of complex numbers.
In a neighborhood $U$ of $a_n$, we can represent $u^*$ as 
\begin{equation} \label{eq:representation}
u^*(x') = u_0(x') e^{i\phi(x')} \qquad \textrm{ for } x'\in U.
\end{equation}
The function $\phi$ is harmonic and is smooth in $U$ and in particular is smooth in $a_n$ (see  \cite[Corollary~I.2]{BBH} and \cite[Formula 2, Chapter VII]{BBH}). 
Upon choosing $\sigma$ small enough we can suppose that $U$ contains the annulus $B_{2\sigma}(a_n)\setminus B_{\sigma}(a_n)$. 
Therefore, the representation~\eqref{eq:representation} holds, in particular, on $\partial B_{\sigma}(a_n)$. 
The smoothness of $\phi$ in $U$ implies that there exists a constant $C$ such that 
\[
\abs{\nabla \phi}\le C \qquad \textrm{ in } U.
\]
Thus, representing~$\phi$ in polar coordinates $(\rho, \theta)$ we have that 
\[
\abs{\partial_\theta \phi_\sharp}\le C\sigma \qquad \textrm{ on } \partial B_{\sigma}(a_n),
\]
where we denote $f_{\sharp}(\rho, \theta) = f(\rho\cos\theta,\rho\sin\theta)$
for any function~$f$.

Now we construct $w_\sigma\colon B_{2\sigma}(a_n)\setminus B_{\sigma}(a_n)\to \mathbb{S}^1$ such that 
\begin{equation}
\label{eq:interpol}
\begin{split}
w_\sigma = u^* \qquad &\textrm{on } \partial B_{2\sigma}(a_n) \\
w_{\sigma} = \alpha_\sigma u_0\qquad &\textrm{on } \partial B_{\sigma}(a_n),
\end{split}
\end{equation}
for some constant~$\alpha_\sigma\in \mathbb{C}$ with~$\abs{\alpha_\sigma}=1$, and
\begin{equation}
\label{eq:interpol-bis}
\frac{1}{2}\int_{B_{2\sigma}(a_n)\setminus B_{\sigma}(a_n)}\abs{D w_\sigma}^2 \d x' \le \pi\log 2 + \mathrm{o}_{\sigma \to 0}(1).
\end{equation}
We let $\bar \phi$ be the mean of $\phi$ over $\partial B_\sigma(a_n)$.
The map $w_\sigma$ is defined as follows (see \cite[Subsection 9.2]{JerrardIgnat_full} for an analogous construction):
\[
w_\sigma(x') = u_{0}(x') e^{i \left[ \bar{\phi} + \left(\frac{\abs{x'-a_n}}{\sigma}-1\right)(\phi_{|\partial B_\sigma(a_n)}-\bar{\phi})\right]}\qquad
\textrm{for } x'\in B_{2\sigma}(a_n)\setminus B_{\sigma}(a_n).
% w_\sigma(x') := u_{0}(x') \exp\left\{i \left[ \bar{\phi} + \left(\frac{\abs{x'-a_n}}{\sigma}-1\right)\left(\phi\left(\frac{\sigma x'}{\abs{x'}}\right)-\bar{\phi}\right)\right]\right\}\qquad \textrm{for } x'\in B_{2\sigma}(a_n)\setminus B_{\sigma}(a_n).
\]
% In particular,
It is easy to show that~$w_\sigma$ satisfies~\eqref{eq:interpol} with $\alpha_\sigma=e^{i\bar{\phi}}$. 
Now we prove \eqref{eq:interpol-bis}. 
First of all, we represent $w_\sigma$ in polar coordinates. We have
\[
(w_{\sigma})_{\sharp}(\rho, \theta) = e^{i\left[\theta + \bar{\phi} + \left(\frac{\rho}{\sigma}-1\right)(\tilde{\phi}_{\sharp}-\bar{\phi})\right]},
\]
where we have set $\tilde{\phi}:= \phi_{|\partial B_\sigma(a_n)}$.
Then, we recall that for a smooth function $f$ there holds
\[
\abs{D f}^2 = \abs{\partial_\rho f_{\sharp}}^2 + \frac{1}{\rho^2}\abs{\partial_\theta f_{\sharp}}^2.
\] 
Therefore (recalling that $\abs{\partial_\theta\tilde{\phi}_{\sharp}}\lesssim \sigma$ on $\partial B_\sigma(a_n)$)
\[
 \begin{split}
 \frac{1}{2}\int_{B_{2\sigma}(a_n)\setminus B_{\sigma}(a_n)}\abs{D w_\sigma}^2 \d x'  &= 
 \frac{1}{2}\int_{0}^{2\pi}\int_{\sigma}^{2\sigma}\left\{\frac{1}{\sigma^2}\abs{\tilde\phi_{\sharp} -\bar{\phi}}^2
 + \frac{1}{\rho^2}\abs{1+ \left(\frac{\rho}{\sigma}-1\right)\partial_{\theta}\tilde{\phi}_{\sharp}}^2\right\}\d \rho\d \theta\\
 &\le C\int_{0}^{2\pi}\abs{\partial_\theta \tilde{\phi}_{\sharp}}^2 \d \theta + \pi \int_{\sigma}^{2\sigma}\frac{1}{\rho}\d \rho + \mathrm{O}_{\sigma	\to 0}(\sigma)\\
 & \le \pi \log 2 + \mathrm{o}_{\sigma \to 0}(1).
 \end{split}
\]
 \medskip
%The details of the construction are in \cite[Subsection 9.2]{JerrardIgnat_full}.
%(Essentially, the construction of~$w_\eps$ is based on 
%interpolating linearly, with respect to the radial
%variable in the annulus, between~$\phi$ and a constant~$\bar{\phi}$,
%defined as the average of~$\phi$ on~$\partial  B_{\sqrt{\eps}}(a_n)$.
%In particular, $\alpha_\eps := \exp(i\bar{\phi})$.)
% based on a linear interpolation in the lifting of $u^*= u_0 e^{i\phi}$ between the constant rotation $\alpha_\eps = e^{i\bar{\phi}}$ ($\bar{\phi}$ denotes the mean of $\phi$ on $\partial B_{\sqrt{\eps}}(a_n)$) and the rotation $e^{i\phi}$.
Given such a~$w_\sigma$, we define
\begin{equation}
W_\sigma(x) := (w_\sigma(x'), 0)\qquad 
\textrm{for }  x=(x',x_3)\in C_{\sqrt{\eps}}(a_n)\setminus C_{\sqrt{\eps}/2}(a_n),
\end{equation}
so as to obtain 
\begin{equation}
\label{eq:energyW}
F_\sigma\left(W_\sigma;C_{2\sigma}(a_n)\setminus C_{\sigma}(a_n)\right) = \frac{1}{2}\int_{B_{2\sigma}(a_n)\setminus B_{\sigma}(a_n)}\abs{D w_\sigma}^2 \d x' \le \pi \log 2 + \mathrm{o}_{\sigma\to 0}(1)
\end{equation}
The sequence $U_{\eps,\sigma}$ is defined to be equal to %the corresponding 
$W_\sigma$ in each~$C_{2\sigma}(a_n)\setminus C_{\sigma}(a_n)$.

Then we work in $C_{\sigma}(a_n)$. Here we consider a map~$\tilde{V}_{\eps,\sigma}\colon C_{\sigma}(a_n)\to\S^2$
such that
\[
F_{\eps}\left(\tilde{V}_{\eps,\sigma};C_{\sigma}(a_n)\right)  = \gamma\left(\sigma,\eps\right)
\]
(i.e., a minimizer for the problem defined in~\eqref{eq:core}),
we consider the rigid motion $R_{\sigma}\colon\R^3\to\R^3$
given by $R_{\sigma}(x):= (\alpha_{\sigma} x^\prime, x_3)$
for~$x= (x^\prime, x_3)\in\R^2\times\R\simeq\C\times\R$
(where~$\alpha_{\sigma}$ is given by~\eqref{eq:interpol}),
% rigid motion $R_\eps$ the rigid motion $R_\eps:= (\alpha_\eps,\textrm{Id})$ (the constant rotation $\alpha_\eps$ defined in \eqref{eq:interpol} in each plane $x_3=t$ for any $t\in (0,1)$ and the identity in the $\hat{e}_3$ direction)
and we define
\begin{equation}
\label{eq:recovery_core}
V_{\eps,\sigma} := R_{\sigma} \tilde{V}_{\eps,\sigma} \qquad \textrm{in } C_{\sigma}(a_n) 
% \begin{cases}
%  R_\eps \tilde{V}_\eps \qquad &\textrm{in } B_{\sqrt{\eps}/2}(a_n) 
%   \qquad \textrm{if } d_n = 1, \\
%  R_\eps \overline{\tilde{V}_\eps} \qquad &\textrm{in } B_{\sqrt{\eps}/2}(a_n) 
%   \qquad \textrm{if } d_n = -1, 
% \end{cases}
\end{equation}
Since the energy~$F_\eps$ is invariant under the rigid transformation~$R_{\sigma}$,
Lemma~\ref{lemma:core} implies that,
along the subsequence~$\eps_k\to 0$ we chose 
at the beginning of the proof, there holds
\begin{equation}
\label{eq:energy_core}
F_{\eps_k}\left(V_{\eps_k,\sigma};C_{\sigma}(a_n)\right)=F_{\eps_k}\left(\tilde{V}_{\eps_k,\sigma};C_{\sigma}(a_n)\right)=\gamma\left(\sigma,\eps_k\right) = \pi\abs{\log\eps_k} -\pi \abs{\log \sigma} + \gamma 
+ \mathrm{r}(\eps_k,\sigma),
\end{equation}
where $\mathrm{r}(\eps_k,\sigma)$ has the property that $\lim_{\sigma\to 0}\lim_{k	\to +\infty}\mathrm{r}(\sigma,\eps_k) = 0$. 
We define~$U_{\eps,\sigma}$ to be equal to~$V_{\eps,\sigma}$
inside each~$C_{\sigma}(a_n)$.
%{\BBB Tecnicamente, tutta la dimostrazione \`e stata fatta
%lungo la sottosuccessione~$\eps_k$ e qui non stiamo pi\`u
%estraendo nessun'altra sottosuccesione: c'\`e scritto~$\eps$
%anzich\'e~$\eps_k$, ma \`e solo per semplicit\`a di notazione.}

To sum up, we have constructed the sequence defined by 
\begin{equation}
\label{eq:recovering}
U_{\eps,\sigma} := 
\begin{cases}
U^* \qquad &\textrm{in } Q_\eps\\
W_\sigma \qquad &\textrm{in } C_{2\sigma}(a_n)\setminus C_{\sigma}(a_n) \quad \textrm{ for }n=1,\ldots,N\\
V_{\eps,\sigma} \qquad &\textrm{in }  C_{\sigma}(a_n)\quad \textrm{ for } n=1,\ldots,N.
\end{cases}
\end{equation}
Along the same subsequence~$\eps_k\to 0$ we have chosen 
previously, we have that (see \eqref{eq:energyU}, \eqref{eq:energyW} and \eqref{eq:energy_core}) 
\begin{equation}
\label{eq:energy_reco}
\begin{split}
F_{\eps_k}(U_{\eps_k,\sigma}) &= F_{\eps}(U^*;Q_\sigma) +\sum_{n=1}^N F_{\eps_k}\left(W_\sigma; C_{2\sigma}(a_n)\setminus C_{\sigma}(a_n)\right) + \sum_{n=1}^{N}F_\eps\left(V_{\eps_k,\sigma}; C_{\sigma}(a_n)\right)\\
&\le  \pi N\abs{\log(\sigma)} + W_g(\a,\db) + N \pi\abs{\log\eps_k} -\pi N\abs{\log \sigma} + N\gamma 
+ \mathrm{r}(\eps_k,\sigma)\\
&\le N \pi\abs{\log\eps_k} + W_g(\a,\db) + N\gamma +\mathrm{r}(\eps_k,\sigma),
\end{split}
\end{equation}
and thus
\begin{equation}
\label{eq:limsup_reco}
\limsup_{k\to +\infty}\left(F_{\eps_k}(U_{\eps_k,\sigma})-N \pi\abs{\log\eps_k}\right)\le W_g(\a,\db) + N\gamma + \mathrm{o}_{\sigma\to 0}(1).
\end{equation}
% Moreover, thanks to 
% ~\eqref{eq:energyW}, \eqref{eq:energy_core} 
% and the H\"older inequality, 
% we have that, for any $p\in [1,2)$ ({\BBB add some detail?}), 
% \begin{equation}
% \label{eq:reco2}
% \limsup_{k\to 0}\norm{U_{\eps_k,\sigma}-U^*}_{W^{1,p}(Q)}\le \mathrm{o}_{\sigma\to 0}(1).
% \end{equation}
By a standard diagonal argument (see Lemma~\ref{lemma:recovery} 
below for details), we find a sequence of positve numbers~$\sigma_k \to 0$
such that~$U_{\eps_k} := U_{\eps_k, \sigma_k}$ satisfies
\begin{equation} \label{eq:limsupenergy}
\limsup_{k\to+\infty} \left(F_{\eps_k}(U_{\eps_k})- N\pi\abs{\log\eps_k}\right)
\le  W_g(\a, \db) + N\gamma.
\end{equation}
It only remains to check that
\begin{equation} \label{eq:limsupfinal}
 J\bar{u}_{\eps_k} \to 2\pi \sum_{n=1}^N \delta_{a_n}
 \qquad \textrm{in } (W^{1,\infty}(\Omega))^\prime,
\end{equation}
where $\bar{u}_{\eps_k}:= \int_0^1 \Pi(U_{\eps_k}) \, \d x_3$.
By and~\eqref{eq:limsupenergy} and Theorem~\ref{th:compactness}, 
we can extract a (non-relabelled) subsequence in such a way that
$J\bar{u}_{\eps_k} \to \mu$ in~$(W^{1,\infty}(\Omega))^\prime$,
where~$\mu$ is a finite sum of Dirac's deltas with integer multiplicities.
The support of~$\mu$ must be contained in the set of points~$\{a_1, \ldots, a_n\}$, because $\bar{u}_{\eps_k}$ coincides
with the canonical harmonic map~$u^*$
in~$\Omega\setminus\cup_{n=1}^N B_{\sigma_k}(a_n)$
(by construction), and hence $J\bar{u}_{\eps_k} = Ju^* = 0$
in~$\Omega\setminus\cup_{n=1}^N B_{\sigma_k}(a_n)$.
The multiplicities of~$\mu$ can be uniquely identified 
by applying the property~\eqref{eq:degree_Jac}, and it turns out
that $\mu = 2\pi \sum_{n=1}^N \delta_{a_n}$,
because the canonical harmonic field~$u^*$ has
a singularity of degree~$d_n = 1$ at each point~$a_n$.
As the limit~$\mu$ is uniquely identified, we deduce that
$J\bar{u}_{\eps_k} \to \mu$ not only along a subsequence,
but for the original sequence as well. 
Therefore, \eqref{eq:limsupfinal} is proved.

% Therefore, using Lemma \ref{lemma:recovery} we conclude that we can extract a further subsequence
% $U_{\eps_k}:= U_{\eps_k,\sigma_k}$ such that 
% \[
% \limsup_{k\to+\infty} \left(F_{\eps_k}(U_{\eps_k})- N\pi\abs{\log\eps_k}\right)
% \le  W_g(\a, \db) + N\gamma,
% \]
% and 
% \[
% U_{\eps_k}\xrightarrow{k\to +\infty}U^* \qquad \textrm{strongly in } \quad W^{1,p}(Q; \R^3)\quad \forall p\in [1,2).
% \]
% This last convergence implies that 
% \[
% \bar{u}_{\eps_k}:= \int_0^1 \Pi(U_\eps) \, \d x_3 \xrightarrow{k\to +\infty}u^* \qquad \textrm{ strongly in } W^{1,p}(\Omega; \R^2)\quad \forall p\in [1,2).
% \]
% Thus, recalling ~\eqref{eq:prejaco} and~\eqref{eq:jaco}, we have that~$j(\bar{u}_\eps)\to j(u^*)$ strongly in~$L^1(\Omega;\R^2)$
% and, hence, $J \bar{u}_\eps \to J u^* = 2\pi \sum_{n=1}^N \delta_{a_n}$
% strongly in~$(W^{1,\infty}(\Omega))^\prime$.

%it is not difficult to
%check that~$U_\eps \to U^*$ strongly in~$W^{1,p}(Q; \R^3)$
%for any~$p\in (1, \, 2)$, and hence
%$\bar{u}_\eps := \int_0^1 \Pi(U_\eps) \, \d x_3 \to u^*$
%strongly in~$W^{1,p}(\Omega;\R^2)$ for any~$p\in (1, \, 2)$.
%Then, using the definition of Jacobian
%(see~\eqref{eq:prejaco} and~\eqref{eq:jaco}),
%we see that~$j(\bar{u}_\eps)\to j(u^*)$ strongly in~$L^1(\Omega;\R^2)$
%and, hence, $J \bar{u}_\eps \to J u^* = 2\pi \sum_{n=1}^N \delta_{a_n}$
%strongly in~$(W^{1,\infty}(\Omega))^\prime$.

So far, we have assumed that~$d_n = 1$ for any~$n$.
However, the proof carries over to the general case,
with minor modifications only.
More precisely, if~$n$ is an index such that~$d_n = -1$,
we define~$U_{\eps,\sigma} := \overline{V_{\eps,\sigma}}$
in~$C_{\sigma}(a_n)$, where~$V_{\eps,\sigma} = (V_{\eps,\sigma,1}, 
V_{\eps,\sigma,2}, V_{\eps,\sigma,3})$
is given by~\eqref{eq:recovery_core} and~$\overline{V_{\eps,\sigma}}$
is defined component-wise as
$\overline{V_{\eps,\sigma}} = (V_{\eps,\sigma,1}, 
-V_{\eps,\sigma,2}, V_{\eps,\sigma,3})$. We modify~$U_{\eps,\sigma}$ 
on~$C_{2\sigma}(a_n)\setminus C_{\sigma}(a_n)$ accordingly.
Effectively, we change the sign of one component of~$U_{\eps,\sigma}$, 
where needed. This process does not affect the energy, but modifies 
the local degree of~$U_{\eps,\sigma}$ so that it agrees with~$d_n$.
\qedhere
\end{step}
\end{proof}

In the last part of the proof of the upper bound we have used the following technical lemma. 

\begin{lemma}[A diagonal argument for recovery sequences]
\label{lemma:recovery}
 Let~$\mathrm{r}\colon \N\times (0, \, +\infty)\to\R$
 be a non-negative function such that
 \begin{equation} \label{almostrecovery}
  \limsup_{\sigma\to 0}\limsup_{k\to+\infty} \mathrm{r}(k, \, \sigma) = 0.
 \end{equation}
 Then, there exists a positive sequence~$\{\sigma_k\}_{k\in\N}$
 such that 
 \[
  \lim_{k\to 0} \sigma_k = 0, \qquad
  \lim_{k\to 0} \mathrm{r}(k, \sigma_k) = 0.
 \]
\end{lemma}
\begin{proof}
 Due to~\eqref{almostrecovery}, there exists a 
 sequence~$\bar{\sigma}_j\to 0$ such that
 \begin{equation} \label{almostrecovery0}
  \limsup_{k\to+\infty} \mathrm{r}(k, \, \bar{\sigma}_j) \leq \frac{1}{j}
  \qquad \textrm{for any } j\in\N.
 \end{equation}
 Let us take~$j=1$ first. By 
 definition of limit, \eqref{almostrecovery0} 
 implies that there exists a strictly increasing
 sequence of indices~$\{\ell(n, \, 1)\}_{n\in\N}$ such that
 \begin{equation*} %\label{almostrecovery1}
  \mathrm{r}(k, \, \bar{\sigma}_1) \leq \frac{1}{n} + 1 \qquad
  \textrm{for any } n\in\N \textrm{ and } 
  k\geq \ell(n, \, 1).
 \end{equation*}
 In a similar way, for~$j=2$ we find a (strictly increasing) 
 subsequence~$\{\ell(n, \, 2)\}_{n\in\N}$
 of~$\{\ell(n, \, 1)\}_{n\in\N}$ such that
 \begin{equation*} %\label{almostrecovery2}
  \mathrm{r}(k, \, \bar{\sigma}_2) \leq \frac{1}{n} + \frac{1}{2} 
  \qquad \textrm{for any } n\in\N \textrm{ and } 
  k\geq \ell(n, \, 2).
 \end{equation*}
 By iterating this argument, we can construct a countable
 family of strictly increasing sequences, $\{\ell(n, \, j)\}_{n\in\N}$
 for any~$j\in\N$, such that each~$\{\ell(n, \, j+1)\}_{n\in\N}$
 is a subsequence of~$\{\ell(n, \, j)\}_{n\in\N}$ and
  \begin{equation} \label{almostrecoveryj}
  \mathrm{r}(k, \, \bar{\sigma}_j) \leq \frac{1}{n} + \frac{1}{j} 
  \qquad \textrm{for any } n\in\N, \ j\in\N \ \textrm{ and } 
  k\geq \ell(n, \, j).
 \end{equation}
 Now, for any~$k\in\N$ with~$k\geq \ell(1,\, 1)$,
 let~$j = j(k)$ be the largest natural number
 such that~$k\geq \ell(j, \, j)$. %By definition,
 The sequence~$\{j(k)\}_{k\in\N}$ is nondecreasing 
 and satisfies $k < \ell(j(k) + 1, \,  j(k) + 1)$, 
 so we must have~$j(k)\to+\infty$ as~$k\to+\infty$.
 We define a sequence~$\{\sigma_k\}_{k\in\N}$ by
 \[
  \sigma_{k} := \bar{\sigma}_{j(k)} 
   \qquad \textrm{if } k\geq \ell(1, \, 1)
 \]
 and~$\sigma_k := \bar{\sigma}_1$ otherwise. Then, $\sigma_k\to 0$
 as~$k\to+\infty$ and
 \[
  \mathrm{r}(k, \, \sigma_{k}) \leq \frac{2}{j(k)} \to 0
  \qquad \textrm{as } k\to+\infty,
 \]
 thanks to~\eqref{almostrecoveryj}.
\end{proof}

\section{Equilibrium Equations}
\label{sec:EqEq}

Our main result, Theorem~\ref{th:minimizers}
in the introduction, is essentially a corollary of
Theorem~\ref{th:compactness}, Remark~\ref{rk:compactness}
and Theorem~\ref{th:Gamma}, combined. However,
we need to show that the limit of a sequence of minimizers 
of~$F_\eps$ is a canonical harmonic map, in the sense
of Definition~\ref{def:chm}. While this result could be
achieved by energy methods, here we adopt a different strategy.
% based on partial differential equations.
More precisely, in this section we derive the equilibrium 
equations for the energy \eqref{eq:energyQ}
and we exploit them to complete the proof of Theorem~\ref{th:minimizers}. 

We let $U\in \mathcal{A}_{G}$ and we take a smooth $\Phi:Q\to \R^3$ such that
$\Phi_{|\partial \Omega\times (0,1)} = 0$. Then, for $t\in (-1,1)$ we consider the map 
\[
 V(t):= \frac{U + t\Phi}{\abs{U+ t \Phi}}.
\]
Note that $V(t)$ is well defined as $\abs{U+ t \Phi}>0$ for sufficiently small $t$ and that $V\in \mathcal{A}_{G}$. 
Moreover $V(0) = U$ and $t\mapsto V(t)$ is differentiable with 
\[
V'(0) = \Phi - \left(U, \Phi\right)U.
\]
We let 
\begin{equation}
\label{eq:f(t)}
f(t):= \frac{1}{2}\int_{Q}\abs{\nabla_\eps V(t)}^2 \d x + \frac{1}{2\eps^2}\int_{\Omega\times \left\{0,1\right\}}\abs{(V(t),\nu)}^2\d x'=: B(t) + \Gamma(t). 
\end{equation} 
On the one hand we have that (see, e.g., \cite[Chapter 1]{LinWang})
\begin{equation}
\label{eq:bulk}
B'(0) = \int_{Q}\nabla_\eps U: \nabla_\eps \Phi \, \d x -\int_{Q}\abs{\nabla_\eps U}^2\left(U,\Phi\right)\d x.
\end{equation}
On the other hand, 
\begin{equation}
\label{eq:Gamma}
\Gamma'(t) = \frac{1}{\eps^2}\int_{\Omega\times \left\{0,1\right\}}\left(V(t),\nu\right)\left(V'(t),\nu\right)\d x',
\end{equation}
and thus
\[
\Gamma'(0) = \frac{1}{\eps^2}\int_{\Omega\times \left\{0,1\right\}}\left(U,\nu\right)\left(\Phi-\left(U,\Phi\right)U, \nu\right)\d x'.
\]
By direct calculation we have that 
\[
\begin{split}
\left(U,\nu\right)\left(\Phi-\left(U,\Phi\right)U, \nu\right) &= \left(\left(U,\nu\right)\Phi,\nu\right)-\left(\left(U,\nu\right)\left(U,\Phi\right)U,\nu\right)\\
& = \left(\Phi, (U,\nu)\nu- (U,\nu)^2U\right) = \left(U,\nu\right)\left(\Phi,\nu-(U,\nu)U\right).
\end{split}
\]
Therefore we have that 
\begin{equation}
\label{eq:Gamma0}
\Gamma'(0) = \frac{1}{\eps^2}\int_{\Omega\times \left\{0,1\right\}}
 \left(U,\nu\right)\left(\Phi,\nu-(U,\nu)U\right)\d x'.
\end{equation}

As a result, the equilibrium condition $B'(0) + \Gamma'(0) = 0$ rewrites as
\begin{equation}
\label{eq:weakEL}
\int_{Q}\nabla_\eps U: \nabla_\eps \Phi \, \d x -\int_{Q}\abs{\nabla_\eps U}^2\left(U,\Phi\right)\d x +
\frac{1}{\eps^2}\int_{\Omega\times \left\{0,1\right\}}
 \left(U,\nu\right)\left(\Phi,\nu-(U,\nu)U\right)\d x' = 0, 
\end{equation}
for any test function~$\Phi\in H^1(Q, \, \R^3)\cap L^\infty(Q, \, \R^3)$
with trace equal to zero on~$\partial\Omega\times (0, \, 1)$
(the regularity assumptions on~$\Phi$ can be relaxed, 
by an approximation argument).

In particular, if we consider $\Phi\in C^{\infty}_c\left(Q;\R^3\right)$ we obtain $B'(0) = 0$, 
%If we impose the equilibrium condition
%\begin{equation}
%\label{eq:equilibrium_cond}
%B'(0) + \Gamma'(0) = 0
%\end{equation}
%and consider $\Phi\in C^{\infty}_c\left(Q;\R^3\right)$, we obtain
%\[
%B'(0) = 0,
%\]
that means 
\begin{equation}
\label{eq:interior_weak}
\int_{Q}\nabla_\eps U :\nabla_\eps\Phi \, \d x = \int_{Q}\abs{\nabla_\eps U}^2 (U,\Phi) \, \d x,
\end{equation}
or, in the distributional sense, 
\[
-\Delta_\eps U = \abs{\nabla_\eps U}^2 U \qquad \textrm{ in } Q, 
\]
where 
\[
\Delta_\eps U := \sum_{k=1}^2 \frac{\partial^2 U}{\partial x_k^2} + \frac{1}{\eta^2(\eps)}\frac{\partial^2 U}{\partial x_3^2}
\]
The last term in \eqref{eq:weakEL} is formally related to a Neumann boundary condition on $\Omega\times \left\{0,1\right\}$. 
More precisely, by a formal (see Remark \ref{rem:regularity} on the regularity issue) integration by parts in \eqref{eq:weakEL}
we obtain that for any
%
%
%
%
%
%In particular, the distribution $-\Delta_\eps U - \abs{\nabla_\eps U}^2 U\in L^2(Q;\R^3)$. 
%Thanks to this information, we obtain that, for any
$\Phi\colon Q\to \R^3$ with $\Phi_{|\partial\Omega\times (0,1)}=0$, 
\[
B'(0) = \int_{\Omega\times \left\{0,1\right\}}\left(\frac{\partial U}{\partial \nu},\Phi\right) \d x'. 
\]
Therefore, from the equilibrium condition $B'(0) +\Gamma'(0)=0$, we obtain that 
\[
\Gamma'(0) = -\int_{\Omega\times \left\{0,1\right\}}\left(\frac{\partial U}{\partial \nu},\Phi\right) \d x'.
\]
Summing up, we interpret the variational equation \eqref{eq:weakEL} as the weak form of the following boundary value problem
\begin{equation}
\label{eq:EL}
\begin{cases}
-\Delta_\eps U = \abs{\nabla_\eps U}^2 U &\qquad \textrm{ in } Q\\
U = G &\qquad \textrm{ on } \partial\Omega\times (0,1)\\
\dfrac{\partial U}{\partial \nu} = -\dfrac{1}{\eps^2}\left(U,\nu\right)\left(\nu-(U,\nu)U\right) &\qquad \textrm{ on } \Omega\times \left\{0,1\right\}.
\end{cases}
\end{equation}
%in the sense of distributions. Equivalently, $U$ satisfies
%\begin{equation}
%\label{eq:weakEL}
% \int_{Q}\nabla_\eps U :\nabla_\eps\Phi \, \d x - \int_{Q}\abs{\nabla_\eps U}^2 (U,\Phi) \, \d x +  \frac{1}{\eps^2}\int_{\Omega\times \left\{0,1\right\}}
% \left(U,\nu\right)\left(\Phi,\nu-(U,\nu)U\right)\d x' = 0 
%\end{equation}
%for any test function~$\Phi\in H^1(Q, \, \R^3)\cap L^\infty(Q, \, \R^3)$
%with trace equal to zero on~$\partial\Omega\times (0, \, 1)$
%(the regularity assumptions on~$\Phi$ can be relaxed, 
%by an approximation argument). 

\begin{remark}
\label{rem:regularity}
It might be worth recalling that weak solutions~$U\in H^1(Q; \S^2)$ 
of~\eqref{eq:EL} cannot be expected to be smooth, in general.
Indeed, Rivi\`ere~\cite{Riviere-harmonic_mostro} showed that,
on a three-dimensional domain, 
there exist weak solutions to the equation $-\Delta U = \abs{\nabla U}^2 U$
(where~$\Delta$ is the standard, isotropic, Laplacian)
that are discontinuous on a dense set.
Minimizers, however, enjoy much better regularity properties.
In particular, %up to a change of variables, we can apply
Schoen and Uhlenbeck's regularity result~\cite{SU}
(see also~\cite{BrezisCoronLieb})
implies that minimizers of~$F_\eps$ are smooth %in~$Q$
except for a locally finite number of singular points.
There are other classes of solutions of~\eqref{eq:EL}
for which a partial regularity theory exist --- 
in particular, the so-called `stationary harmonic' maps,
see~\cite{Evans-PartialReg, Bethuel-PartialReg}.
The interested reader is referred, e.g., to~\cite{Moser}
for a review of this topic.
% partial regularity results for harmonic maps.
\end{remark}

 Concerning the regularity of minimizers of~$F_\eps$,
we have the following result. 

\begin{prop}\label{prop:smooth}
 Assume that the boundary datum~$G$ satisfies~\eqref{eq:bddatum},
 for some~$g\in H^{1/2}(\partial\Omega;\S^1)$.
 Then, for any~$\eps> 0$, minimizers~$U_\eps^*$ of~$F_\eps$ in~$\mathcal{A}_G$
 are smooth in the interior of~$Q$
 and~$U_{\eps,3}^* := (U_\eps^*, \hat{e}_3)$ has constant 
 sign in~$Q$.
\end{prop}
\begin{proof}
 A regularity result by Giaquinta and
 Sou\oldv{c}ek~\cite{GiaquintaSoucek} implies that
 minimizers of~$F_\eps$ that take values in the half-sphere
 $\S^2_+:=\{y\in\S^2\colon y_3\geq 0\}$ are smooth in the 
 interior of~$Q$. Now, let~$U_\eps^*$ be any minimizer of~$F_\eps$
 in the class~$\mathcal{A}_G$. The 
 map~$V_\eps :=(U^*_{\eps,1}, U^*_{\eps,2}, \abs{U^*_{\eps,3}})$
 belongs to~$\mathcal{A}_G$ and is still a minimizer of~$F_\eps$.
 Therefore, $V_\eps$ is smooth (because of~\cite{GiaquintaSoucek})
 and satisfies the system~\eqref{eq:EL}. In particular,
 $V_{\eps,3} := (V_\eps, \hat{e}_3) = \abs{U^*_{\eps,3}}$ 
 is a smooth, non-negative solution of 
 $-\Delta(V_{\eps,3}) = V_{\eps,3}\abs{\nabla V_\eps}^2\geq 0$.
 By the strong maximum principle, $V_{\eps,3}$ is strictly
 positive inside~$Q$ and the proposition follows.
\end{proof}

Proposition~\ref{prop:smooth} does \emph{not}
guarantee that minimizers of~$F_\eps$ are smooth 
up to the boundary of~$Q$: we cannot apply the results
of~\cite{GiaquintaSoucek} to obtain boundary regularity,
because of the surface term in the functional~$F_\eps$.

We can now give the proof of our main result.
\begin{proof}[Proof of Theorem~\ref{th:minimizers}]
 Let~$U_\eps^*$ be a minimizer of~$F_\eps$
 in the class~$\mathcal{A}_G$. By Theorem~\ref{th:Gamma}
 (in particular, Statement~$(ii)$), we can construct a sequence 
 of comparison maps~$U_\eps\in\mathcal{A}_g$ such that
 \[
  F_\eps(U^*_\eps) \leq F_\eps(U_\eps) \leq d\pi \abs{\log\eps} + C,
 \]
 where~$d$ is the degree of the boundary datum 
 (which is assumed to be positive) and~$C$ is 
 an~$\eps$-independent constant. Then, Theorem~\ref{th:compactness}
 and Remark~\ref{rk:compactness} imply that,
 up to extraction of a non-relabelled subsequence,
 $U^*_\eps$ converges %(strongly in~$L^q(Q)$ for any finite~$q$)
 to a map of the form~$U^*(x) = (u^*(x^\prime), 0)$;
 the convergence is strong in~$L^q(Q)$ for any finite~$q$.
 Moreover, $u^*\in W^{1,p}(\Omega;\S^1)$ for any~$p\in (1, \, 2)$,
 and that~$u^*$ is locally in~$H^1$ away from 
 a finite number of points~$a^*_1, \ldots, a^*_d$,
 which are topological singularities of degree~$1$.
 The $\Gamma$-compactness result, Theorem~\ref{th:Gamma},
 implies that~$a^*_1, \ldots, a^*_d$ minimize
 the Renormalized Energy.
 
 To complete the proof, it only remains to show 
 that~$u^*$ is the canonical harmonic field
 associated with~$\a^*=(a^*_1, \ldots, a^*_d)$ and~$g$
 --- more precisely, it remains to see that~$\div j(u^*) = 0$.
 To this end, we fix a function~$\varphi\in C^\infty_{\mathrm{c}}(\Omega)$
 and we test the equation~\eqref{eq:weakEL} against~$\Psi(x) := (-\varphi(x^\prime) U_2(x), \varphi(x^\prime) U_1(x), \, 0)$.
 This map is orthogonal to both~$U$ and~$e_3$
 at each point. After standard computations, we obtain
 \begin{equation*}
 \sum_{k=1}^2 \int_{Q} \left(-U^*_{\eps,2} \, \partial_k U^*_{\eps,1}
  + U^*_{\eps,1}\,\partial_k U^*_{\eps,2} \right)
  \partial_k \varphi \, \d x  = 0 
 \end{equation*} 
 In terms of~$\bar{u}_\eps^* := \int_0^1\Pi(U_\eps^*)\,\d x_3$,
 this equation rewrites as
 $\div j(\bar{u}_\eps^*) = 0$ in the sense of distributions.
 Due to~\eqref{eq:limitmap}, we can pass to the limit as~$\eps\to 0$
 and conclude that~$\div j(u^*) = 0$. Therefore, $u^*$
 is the canonical harmonic map associated with~$(\a^*, \, g)$.
\end{proof}

\subsection*{Acknowledgements.}
Both authors were supported by GNAMPA-INdAM.
G.C. acknowledges the support of the University of Verona
under the project RIBA 2019 No.~RBVR199YFL ``Geometric
Evolution of Multi Agent Systems'', of the 
\emph{Agence National des la Recherche}
under the projet ANR-22-CE40-0006
``Singularities of energy-minimizing vector-valued maps'',
and of the \emph{Leverhulme Trust} under the Research Project 
ORPG-9787 ``Unravelling the Mysteries of Complex 
Nematic Solution Landscapes''.
A.S. acknowledges the partial support of the MIUR-PRIN Grant 2017 ''Variational methods for stationary and evolution problems with singularities and interfaces''.

\bibliographystyle{plain}
\bibliography{dim_redu}

\begin{thebibliography}{10}

\bibitem{AlamaBronsardGolovaty}
Stan Alama, Lia Bronsard, and Dmitry Golovaty.
\newblock Thin film liquid crystals with oblique anchoring and boojums.
\newblock {\em Annales de l'Institut Henri Poincaré C, Analyse non linéaire},
  37(4):817--853, 2020.

\bibitem{ABO2}
G.~Alberti, S.~Baldo, and G.~Orlandi.
\newblock Variational convergence for functionals of {G}inzburg-{L}andau type.
\newblock {\em Indiana Univ. Math. J.}, 54(5):1411--1472, 2005.

\bibitem{AlicandroPonsiglione}
R.~Alicandro and M.~Ponsiglione.
\newblock Ginzburg-{L}andau functionals and renormalized energy: a revised
  {$\Gamma$}-convergence approach.
\newblock {\em J. Funct. Anal.}, 266(8):4890--4907, 2014.

\bibitem{BallCime}
John~M. Ball.
\newblock Liquid crystals and their defects.
\newblock In {\em Mathematical thermodynamics of complex fluids}, volume 2200
  of {\em Lecture Notes in Math.}, pages 1--46. Springer, Cham, 2017.

\bibitem{BarberoDurand}
{Barbero, G.} and {Durand, G.}
\newblock On the validity of the {R}apini-{P}apoular surface anchoring energy
  form in nematic liquid crystals.
\newblock {\em J. Phys. France}, 47(12):2129--2134, 1986.

\bibitem{Bethuel-PartialReg}
F.~Bethuel.
\newblock On the singular set of stationary harmonic maps.
\newblock {\em manuscripta mathematica}, 78(1):417--443, 1993.

\bibitem{BBH}
F.~Bethuel, H~Brezis, and F~H{\'e}lein.
\newblock {\em Ginzburg-{L}andau vortices}.
\newblock Progress in Nonlinear Differential Equations and their Applications,
  13. Birkh\"auser Boston, Inc., Boston, MA, 1994.

\bibitem{BethuelBrezisOrlandi}
F.~Bethuel, H.~Brezis, and G.~Orlandi.
\newblock Asymptotics for the {G}inzburg-{L}andau equation in arbitrary
  dimensions.
\newblock {\em J. Funct. Anal.}, 186(2):432--520, 2001.

\bibitem{BethuelOrlandiSmets-Annals}
F.~Bethuel, G.~Orlandi, and D.~Smets.
\newblock Convergence of the parabolic ginzburg–landau equation to motion by
  mean curvature.
\newblock {\em Ann. Math.}, 163(1):37--163, 2006.

\bibitem{BoutetdeMonvel}
A.~Boutet~de Monvel-Berthier, V.~Georgescu, and R.~Purice.
\newblock A boundary value problem related to the {G}inzburg-{L}andau model.
\newblock {\em Comm. Math. Phys.}, 142(1):1--23, 1991.

\bibitem{braides-beginner}
Andrea Braides.
\newblock {\em {$\Gamma$}-convergence for beginners}, volume~22 of {\em Oxford
  Lecture Series in Mathematics and its Applications}.
\newblock Oxford University Press, Oxford, 2002.

\bibitem{BrezisCoronLieb}
H.~Brezis, J.-M. Coron, and E.~H. Lieb.
\newblock Harmonic maps with defects.
\newblock {\em Comm. Math. Phys.}, 107(4):649--705, 1986.

\bibitem{brezis}
Haim Brezis.
\newblock {\em Functional analysis, {S}obolev spaces and partial differential
  equations}.
\newblock Universitext. Springer, New York, 2011.

\bibitem{BM21}
Ha\"{\i}m Brezis and Petru Mironescu.
\newblock {\em Sobolev maps to the circle---from the perspective of analysis,
  geometry, and topology}, volume~96 of {\em Progress in Nonlinear Differential
  Equations and their Applications}.
\newblock Birkh\"{a}user/Springer, New York, [2021] \copyright 2021.

\bibitem{AG-Reno}
G.~Canevari and A.~Segatti.
\newblock Dynamics of {G}inzburg-{L}andau vortices for vector fields on
  surfaces.
\newblock Preprint arXiv~https://arxiv.org/abs/2108.01321, 2021.

\bibitem{AGM-VMO}
G.~Canevari, A.~Segatti, and M.~Veneroni.
\newblock Morse's index formula in {VMO} for compact manifolds with boundary.
\newblock {\em J. Funct. Anal.}, 269(10):3043--3082, 2015.

\bibitem{gamma-discreto}
Giacomo Canevari and Antonio Segatti.
\newblock Defects in nematic shells: a {$\Gamma$}-convergence
  discrete-to-continuum approach.
\newblock {\em Arch. Ration. Mech. Anal.}, 229(1):125--186, 2018.

\bibitem{dg}
P.~G. De~Gennes.
\newblock {\em The physics of Liquid Crystals}.
\newblock Clarendon Press, Oxford, Oxford, 1974.

\bibitem{Eri66}
J.~L. Ericksen.
\newblock {Inequalities in Liquid Crystal Theory}.
\newblock {\em The Physics of Fluids}, 9(6):1205--1207, 06 1966.

\bibitem{ericksen76}
J.L. Ericksen.
\newblock Equilibrium theory of liquid crystals.
\newblock volume~2 of {\em Advances in Liquid Crystals}, pages 233--298.
  Elsevier, 1976.

\bibitem{Evans-PartialReg}
L.~C. Evans.
\newblock Partial regularity for stationary harmonic maps into spheres.
\newblock {\em Arch. Rational Mech. Anal.}, 116(2):101--113, 1991.

\bibitem{frank58}
F.~C. Frank.
\newblock I. liquid crystals. on the theory of liquid crystals.
\newblock {\em Discuss. Faraday Soc.}, 25:19--28, 1958.

\bibitem{GiaquintaSoucek}
M.~Giaquinta and J.~Sou\oldv{c}ek.
\newblock Harmonic maps into a hemisphere.
\newblock {\em Annali della Scuola Normale Superiore di Pisa - Classe di
  Scienze}, 12(1):81--90, 1985.

\bibitem{GMS15}
Dmitry Golovaty, Jos\'{e}~Alberto Montero, and Peter Sternberg.
\newblock Dimension reduction for the {L}andau-de {G}ennes model in planar
  nematic thin films.
\newblock {\em J. Nonlinear Sci.}, 25(6):1431--1451, 2015.

\bibitem{GMS17}
Dmitry Golovaty, Jos\'{e}~Alberto Montero, and Peter Sternberg.
\newblock Dimension reduction for the {L}andau--de {G}ennes model on curved
  nematic thin films.
\newblock {\em J. Nonlinear Sci.}, 27(6):1905--1932, 2017.

\bibitem{HKL86}
Robert Hardt, David Kinderlehrer, and Fang-Hua Lin.
\newblock Existence and partial regularity of static liquid crystal
  configurations.
\newblock {\em Comm. Math. Phys.}, 105(4):547--570, 1986.

\bibitem{JerrardIgnat_full}
R.~Ignat and R.~L. Jerrard.
\newblock Renormalized energy between vortices in some {G}inzburg-{L}andau
  models on 2-dimensional {R}iemannian manifolds.
\newblock {\em Arch. Ration. Mech. Anal.}, 239(3):1577--1666, 2021.

\bibitem{IgnatKurzke}
Radu Ignat and Matthias Kurzke.
\newblock Global {J}acobian and {$\Gamma$}-convergence in a two-dimensional
  {G}inzburg-{L}andau model for boundary vortices.
\newblock {\em Journal of Functional Analysis}, 280(8):108928, 2021.

\bibitem{Jerrard}
R.~L. Jerrard.
\newblock Lower bounds for generalized {G}inzburg-{L}andau functionals.
\newblock {\em SIAM J. Math. Anal.}, 30(4):721--746, 1999.

\bibitem{JerrardSoner-GL}
R.~L. Jerrard and H.~M. Soner.
\newblock The {J}acobian and the {G}inzburg-{L}andau energy.
\newblock {\em Cal. Var. Partial Differential Equations}, 14(2):151--191, 2002.

\bibitem{LinRiviere}
F.-H. Lin and T.~Rivi{\`e}re.
\newblock Complex {G}inzburg-{L}andau equations in high dimensions and
  codimension two area minimizing currents.
\newblock {\em J. Eur. Math. Soc. (JEMS)}, 1(3):237--311, 1999.

\bibitem{LinWang}
Fanghua Lin and Changyou Wang.
\newblock {\em The analysis of harmonic maps and their heat flows}.
\newblock World Scientific Publishing Co. Pte. Ltd., Hackensack, NJ, 2008.

\bibitem{Moser}
R.~Moser.
\newblock {\em Partial regularity for harmonic maps and related problems}.
\newblock World Scientific Publishing Co. Pte. Ltd., Hackensack, NJ, 2005.

\bibitem{NapVer12L}
G.~Napoli and L.~Vergori.
\newblock Extrinsic curvature effects on nematic shells.
\newblock {\em Phys. Rev. Lett.}, 108(20):207803, 2012.

\bibitem{NapVer12E}
G.~Napoli and L.~Vergori.
\newblock Surface free energies for nematic shells.
\newblock {\em Phys. Rev. E}, 85(6):061701, 2012.

\bibitem{Novack}
Michael~R. Novack.
\newblock Dimension reduction for the {L}andau-de {G}ennes model: The vanishing
  nematic correlation length limit.
\newblock {\em SIAM Journal on Mathematical Analysis}, 50(6):6007--6048, 2018.

\bibitem{oseen33}
C.W. Oseen.
\newblock The theory of liquid crystals.
\newblock {\em Trans. Faraday Soc.}, 29(140):883--899, 1933.

\bibitem{Riviere-harmonic_mostro}
T.~Rivi{\`e}re.
\newblock Everywhere discontinuous harmonic maps into spheres.
\newblock {\em Acta Math.}, 175(2):197--226, 1995.

\bibitem{Sandier}
{\'E}.~Sandier.
\newblock Lower bounds for the energy of unit vector fields and applications.
\newblock {\em J. Funct. Anal.}, 152(2):379--403, 1998.
\newblock see Erratum, ibidem 171, 1 (2000), 233.

\bibitem{SS-book}
{\'E}.~Sandier and S.~Serfaty.
\newblock {\em Vortices in the magnetic {G}inzburg-{L}andau model}.
\newblock Progress in Nonlinear Differential Equations and their Applications,
  70. Birkh\"auser Boston, Inc., Boston, MA, 2007.

\bibitem{SS-GF}
Etienne Sandier and Sylvia Serfaty.
\newblock Gamma-convergence of gradient flows with applications to
  {G}inzburg-{L}andau.
\newblock {\em Comm. Pure Appl. Math.}, 57(12):1627--1672, 2004.

\bibitem{SU}
R.~Schoen and K.~Uhlenbeck.
\newblock Boundary regularity and the {D}irichlet problem for harmonic maps.
\newblock {\em J. Differential Geom.}, 18(2):253--268, 1983.

\bibitem{struwe94}
M~Struwe.
\newblock On the asymptotic behavior of minimizers of the ginzburg-landau model
  in 2 dimensions.
\newblock {\em Diff. and Int. Equations}, 7(6):1631--1624, 1994.

\bibitem{Virga94}
E.~G. Virga.
\newblock {\em Variational theories for liquid crystals}, volume~8 of {\em
  Applied Mathematics and Mathematical Computation}.
\newblock Chapman \& Hall, London, 1994.

\bibitem{zocher33}
H.~Zocher.
\newblock The effect of a magnetic field on the nematic state.
\newblock {\em Trans. Faraday Soc.}, 29:945--957, 1933.

\end{thebibliography}
\end{document}